\documentclass[11pt]{amsart}

\usepackage{amsmath, amssymb, amsfonts, amsthm}
\usepackage{graphicx}
\usepackage{color}
\usepackage{bm}
\usepackage{subfigure}
\usepackage{graphicx}
\usepackage{mathabx}
\usepackage{multirow}
\usepackage{setspace}
\usepackage[title]{appendix}
\usepackage{mathrsfs} 
\usepackage{enumerate}

\usepackage{tikz}
\usepackage{tikz-cd}
\usepackage{wrapfig}
\usepackage{float}

\usetikzlibrary{matrix,arrows,decorations.pathmorphing}

\usepackage{amsthm}
\theoremstyle{plain}
\newtheorem{theorem}{Theorem}[section]

\theoremstyle{plain}
\newtheorem{lemma}[theorem]{Lemma}

\theoremstyle{plain}
\newtheorem{Prop}{Proposition}[section]

\theoremstyle{plain}
\newtheorem{corollary}[theorem]{Corollary}

\theoremstyle{plain}

\theoremstyle{definition}
\newtheorem{Def}{Definition}[section]
\theoremstyle{remark} 
\newtheorem{remark}{Remark}

\theoremstyle{definition}
\newtheorem{notation}{Notation}

\theoremstyle{definition} 
\newtheorem{Example}{Example}

\def\ed{\mathrm{d}}

\def\I{\mathcal I}

\def\J{\mathcal J}
\def\B{\mathcal B}

\def\R{\mathbb R}
\def\C{\mathbb C}
\def\P{\mathbb P}

\def\E{\mathbb E}

\def\e{{\bf e}}

\def\g{\mathfrak g}

\def\so{\mathfrak{so}}

\def\gl{\mathfrak{gl}}
\def\sl{\mathfrak{sl}}

\def\rank{{\rm rank}}

\def\diag{{\rm diag}}

\def\Hom{{\rm Hom}}

\def\W{\wedge}
\def\<{\langle}
\def\>{\rangle}

\def\SO{{\rm SO}}

\def\GL{{\rm GL}}

\def\A{\textbf{A}}
\def\Be{\textbf{B}}

\def\lbb{[\![}
\def\rbb{]\!]}

\def\({\left(}
\def\){\right)}

\def\bs{\boldsymbol}

\DeclareFontFamily{U}{MnSymbolC}{}
\DeclareSymbolFont{MnSyC}{U}{MnSymbolC}{m}{n}
\DeclareFontShape{U}{MnSymbolC}{m}{n}{
    <-6>  MnSymbolC5
   <6-7>  MnSymbolC6
   <7-8>  MnSymbolC7
   <8-9>  MnSymbolC8
   <9-10> MnSymbolC9
  <10-12> MnSymbolC10
  <12->   MnSymbolC12}{}
\DeclareMathSymbol{\lefthook}{\mathbin}{MnSyC}{'270}

\title[Geometry of B\"acklund Transformations I]{Geometry of B\"acklund Transformations I: Generality}
\author{Yuhao Hu} 
\address{Department of Mathematics, 395 UCB, University of
Colorado, Boulder, CO 80309-0395}
\email{Yuhao.Hu@colorado.edu}


\begin{document}

\maketitle

\begin{abstract}
	Using \'Elie Cartan's method of equivalence, we prove an upper bound for the generality of generic rank-1 B\"acklund transformations relating two hyperbolic Monge-Amp\`ere systems. In cases when the B\"acklund transformation admits a symmetry group whose orbits have codimension 1, 2 or 3, we obtain classification results and new examples of auto-B\"acklund transformations. 

\end{abstract}

\setcounter{tocdepth}{1}
\tableofcontents

\section{Introduction}

In 1882, the Swedish mathematician Albert V. B\"acklund proved the result (see \cite{backlund1882}, \cite{BGG} or \cite{chern-terng}): \emph{Given a surface with a constant Gauss curvature $K<0$ in $\E^3$, 
one can construct, by solving ODEs, a 1-parameter family of new surfaces  in $\E^3$ with the Gauss curvature $K$.}
This is the origin of the term ``B\"acklund transformation".

Classically, a \emph{B\"acklund transformation} is a PDE system $\B$ that relates solutions
of two other PDE systems $\mathcal{E}_1$ and $\mathcal{E}_2$. Moreover, such a relation must satisfy the 
property: given a solution $u$ of $\mathcal{E}_1$ (resp. $\mathcal{E}_2$), substituting it in $\B$, one 
would obtain a PDE system
whose solutions can be found by ODE methods and produce solutions of $\mathcal{E}_2$ (resp. $\mathcal{E}_1$). If, in addition, $\mathcal{E}_1$ and $\mathcal{E}_2$
are contact equivalent to each other, then the corresponding B\"acklund
transformation is called an \emph{auto-B\"acklund transformation}.

For example, the \emph{Cauchy-Riemann system}
	\begin{equation}
		\left\{
			\begin{array}{l}
				u_x - v_y = 0,\\[0.8em]
				u_y+ v_x = 0
			\end{array}
		\right.	\label{CR}
	\end{equation}
is an auto-B\"acklund transformation; it relates solutions of the Laplace equation $\Delta z = 0$ for $z(x,y)$ in the following
way: If $u$ satisfies $\Delta u= 0$, then, substituting it in \eqref{CR}, one obtains a compatible first-order system
for $v$, whose solutions can be found by ODE methods and satisfy $\Delta v= 0$, and \emph{vice versa}.

As another example, consider the system of nonlinear equations
	\begin{equation}
		\left\{
			\begin{array}{l}
				z_x - \bar z_x = \lambda \sin(z+\bar z),\\[0.8em]
				z_y+\bar z_y = \lambda^{-1}\sin(z-\bar z),
			\end{array}
		\right.	\label{SGT}
	\end{equation}
where $\lambda$ is a nonzero constant. One can easily verify that \eqref{SGT} is an auto-B\"acklund
transformation relating solutions of the \emph{sine-Gordon equation}
	\begin{equation}
		u_{xy} = \dfrac{1}{2} \sin(2u). \label{SG}
	\end{equation}
The system \eqref{SGT} can be derived from
the classical auto-B\"acklund transformation relating surfaces  in $\E^3$ with a 
negative constant Gauss curvature. For details, see \cite{chern-terng}.

In general, there may seem to be very few restrictions on the types of PDE systems that admit a B\"acklund
transformation. In addition to the elliptic and hyperbolic examples mentioned above, a B\"acklund transformation
may exist relating solutions of a parabolic equation (see \cite{nimmo82}) or an equation of order higher than 2, for example, the KdV equation (see \cite{wahlquist1973}).
Furthermore, two PDE systems being B\"acklund-related need not be contact equivalent
to each other (see \cite{CI09}).

The importance of B\"acklund transformations may, in part, be viewed through their relation to surface geometry
and mathematical physics. On the geometry side, B\"acklund transformations allow one to obtain new surfaces with
prescribed geometric properties from old. For a variety of such examples, see \cite{Rogers}. On the mathematical physics side, a prototypical result is that the B\"acklund transformation \eqref{SGT}, when
applied to the trivial solution $z(x,y) = 0$ of \eqref{SG}, yields a 2-parameter family of 1-soliton solutions of the sine-Gordon equation
(see \cite{terng2000geometry}).
More elaborate techniques have since been developed to find the so-called multi-soliton solutions of nonlinear PDE systems (for example, the KdV equation), using B\"acklund transformations.

An ultimate goal of studying B\"acklund transformations is solving
the \emph{B\"acklund problem}, which was considered by \'E. Goursat in \cite{Goursat1925}: \begin{quote}\emph{Find all pairs of systems of PDEs whose solutions are related by a B\"acklund transformation.} \end{quote} 
Although this problem remains largely unsolved, recent works of Clelland and Ivey (\cite{C01},\cite{CI05} and \cite{CI09}) and those of Anderson and Fels (\cite{AF12}, \cite{AF15})
have pointed out new directions for studying B\"acklund transformations. Instead of aiming at constructing new examples or finding techniques of calculating explicit solutions to PDE systems, they work in a geometric setting that is
natural to the study of structural properties (of B\"acklund transformations) that are invariant under 
contact transformations. Under such settings, a complete classification of B\"acklund transformations, at least in 
certain cases, is possible by using \'E. Cartan's method of equivalence.

The current work is concerned with the geometric aspect of B\"acklund transformations, not so much in the sense of relating to classical surface geometry, as in that of seeing B\"acklund transformations as geometric objects and studying their invariants.

More specifically, we study nontrivial rank-$1$ B\"acklund transformations (see Definition \ref{Backlundrank}) relating a pair of hyperbolic Monge-Amp\`ere systems.
Since many classical examples belong to this category, it is highly desirable to have a complete classification of B\"acklund transformations of this kind. 
In \cite{C01}, by establishing a $G$-structure associated to a B\"acklund transformation, Clelland approached the classification problem using Cartan's
method of equivalence, restricting to the case when all local invariants of the structure are constants (a.k.a the homogeneous case). Her classification found 15 types, within which 11 are analogues of
the classical B\"acklund transformation between surfaces in $\E^3$ with a negative constant Gauss curvature. 

Since homogeneous structures, up to equivalence, depend only on constants, the following question remains to be answered: {\it What kind of initial data do we need
to specify in order to determine a rank-1 B\"acklund transformation relating two hyperbolic Monge-Amp\`ere systems?} 

In Section \ref{generality}, in the generic case, we use the method of equivalence (see \cite{GardnerEquiv}, \cite{Br14}) to prove an upper bound for the magnitude of such initial data:
 \begin{quote}
 \emph{To determine a generic B\"acklund transformation
relating two hyperbolic Monge-Amp\`ere systems, it is sufficient to specify at most $6$ functions of $3$ variables.}
\end{quote}
It is an immediate consequence of our theorem that most hyperbolic Monge-Amp\`ere systems are not related 
to any system of the same type by a generic B\"acklund transformation.

A major difficulty in using the method of equivalence to classify B\"acklund transformations lies in
verifying the compatibility of large systems of polynomial equations, whose variables are the B\"acklund structure invariants and their covariant derivatives. However, we found that such calculation becomes much more manageable when we set two
structure invariants to be specific constants (see Section \ref{rk1example}). The corresponding B\"acklund transformations 
are either homogeneous 
(corresponding to a case already classified in \cite{C01}) or of cohomogeneity 1, 2 or 3. 
In the cohomogeneity-1 case, we obtain an auto-B\"acklund transformation of a homogeneous Euler-Lagrange system
that is contact equivalent to the equation
		\begin{equation*} \tag{\ref{Cohomog1PDE}}
			(A^2-B^2)(z_{xx} - z_{yy}) + 4ABz_{xy} = 0,
		\end{equation*}
where $A = 2z_x+y$ and $B = 2z_y - x$.
The cohomogeneity-2 case has a subcase that arises when
the `free derivatives' associated to the structure are expressed in terms of the primary invariants. 
In this subcase, the corresponding Lie algebra of symmetry must be 
of the form $\mathfrak{q}\oplus \R$ where $\mathfrak{q}$ is either $\sl(3,\R)$, 
$\so(3,\R)$ or the solvable 3-dimensional Lie algebra with basis $\{x_i\}_{i = 1}^{3}$ satisfying
		\[
			[x_2,x_3]  = x_1, \quad [x_3,x_1] = x_2, \quad [x_1, x_2] = 0.
		\]
In particular, when $\mathfrak{q}$ is solvable, we obtain an auto-B\"acklund transformation; the underlying Monge-Amp\`ere system 
is Euler-Lagrange, of cohomogeneity 1, and contact equivalent to the equation
		\begin{equation*}\tag{\ref{Cohomog2PDE}}
			(A^2-B^2)(z_{xx}-z_{yy})+4ABz_{xy} + (A^2+B^2)^2 = 0,
		\end{equation*}
where $A = z_x-y$, $B = z_y + x$.

\section{Definitions and Notations}

In this section, we present some definitions and notations to be used later.

\subsection{Exterior Differential Systems {\rm(c.f. \cite{BCG})}}

   {\Def  Let $M$ be a smooth manifold, $\I\subset \Omega^*(M)$ a graded ideal that is closed under exterior differentiation. The pair $(M,\I)$ is said to be an \emph{exterior differential system} with
		space $M$ and differential ideal $\I$.}\vskip 2mm
		
 Given an exterior differential system $(M,\I)$, 
  we use $\I^k$ to denote the degree-$k$ piece of $\I$, namely, $\I^k = \I\cap \Omega^k(M)$, where $\Omega^k(M)$ stands for the $C^\infty(M)$-module of differential $k$-forms on $M$.
  If the rank of $\I^k$, restricted to each point, is locally a constant, then the elements of $\I^k$ are precisely smooth sections of a vector bundle denoted by $I^k$.
  
  {\Def An \emph{integral manifold} of an exterior differential system $(M,\I)$ is an immersed submanifold $i: N\hookrightarrow M$ satisfying $i^*\phi = 0$ for any $\phi\in \I$.}\vskip 2mm
	
 Intuitively, an exterior differential system is a coordinate-independent way to express a PDE system; an integral manifold, usually with a certain independence condition satisfied,
   corresponds to a solution of the PDE system. 		

  {\Def  Two exterior differential systems $(M,\I)$ and $(N,\J)$ are said to be \emph{equivalent up to diffeomorphism}, or \emph{equivalent}, for brevity, if there exists a diffeomorphism $\phi: M\rightarrow N$ such that 
				$\phi^*\J = \I.$   Such a $\phi$ is called an \emph{equivalence} between the two systems. An equivalence between $(M,\I)$ and itself is called a \emph{symmetry} of $(M,\I)$.

  {\Def Let $\pi:N\rightarrow M$ be a submersion.  
		A $p$-form $\omega\in \Omega^p(N)$ is said to be 
		\emph{$\pi$-semi-basic} if, for any $x\in N$, $\omega|_x\in \pi^*(\Lambda^p(T^*M))$.}		
	
{\Def Let $M$ be a smooth manifold.
		Let $E\subset \Lambda^k(T^*M)$ be a vector subbundle, and $X$ a smooth vector field defined on $M$. 
		We say that $E$ is \emph{invariant under
		the flow of $X$} if, for any $($smooth$)$ local
		section $\omega:U\rightarrow E$, where $U\subset M$ is open, 
		the Lie derivative $\mathcal{L}_X\omega$
		remains a section of $E$ (over $U$).		\label{InvariantBundle}}
	
    {\notation Let $\lbb\theta_1,\ldots,\theta_\ell\rbb$ denote the vector subbundle
							 of $\Lambda^k(T^*U)$ generated 
	by differential forms $\theta_1,\ldots,\theta_\ell$ (defined on $U$, an open subset of a smooth manifold) 
	of the same degree $k$.	}

\subsection{Hyperbolic Monge-Amp\`ere Systems {\rm (c.f. \cite{BGG})}}\	
	
	Among second order PDEs for 1 unknown function of 2 independent variables, \emph{Monge-Amp\`ere equations}
are those of the form
	\begin{equation}
		A(z_{xx}z_{yy} - z_{xy}^2)+Bz_{xx}+2Cz_{xy}+Dz_{yy}+E = 0,	\label{MAequation}
	\end{equation}
where $A,B,C,D,E$ are functions of $x,y,z,z_x,z_y$. A Monge-Amp\`ere equation \eqref{MAequation} is said to be
\emph{elliptic} (resp., \emph{hyperbolic}, \emph{parabolic}) if $AE - BD+ C^2$ is negative (resp.,  positive, zero).

A Monge-Amp\`ere equation can be formulated as an exterior differential system on a contact manifold.
In the hyperbolic case, 
we follow \cite{BGH1} to give the following definition.
				 
{\Def A \emph{hyperbolic Monge-Amp\`ere system} $(M,\I)$ is an exterior differential system, where $M$ is a 5-manifold, 
	$\I$ being locally algebraically generated by $\theta\in \I^1$ and $\ed\theta, \Omega\in \I^2$ satisfying
	\begin{enumerate}[\quad (1)]	
		\item{$\theta\W (\ed\theta)^2 \ne 0$;}
		
		\item{$\lbb \ed\theta, \Omega\rbb$, modulo $\theta$, has rank $2$;}
		
		\item{$(\lambda \ed\theta+ \mu \Omega)^2\equiv 0 \mod \theta$ has two distinct solutions $[\lambda_i:\mu_i]\in \R\P^1$ $(i = 1,2)$.}
	\end{enumerate}	
\label{hyperbolicMAdefinition}}\

Here, condition $(3)$, in particular, characterizes hyperbolicity: Each 
	integral surface of $(M,\I)$ is foliated by two distinct families of characteristics.
	
{\Def Consider a hyperbolic Monge-Amp\`ere system $(M,\I)$. A local coframing $\bs\theta = (\theta^0,\theta^1,\ldots,\theta^4)$ defined on an open neighborhood $U\subset M$ is said to be \emph{$0$-adapted} if, on $U$,
					\[
						\I = \<\theta^0,\theta^1\W\theta^2,\theta^3\W\theta^4\>.
					\]
	\label{MAadmissible}}		
			
The condition for $0$-adaptedness as defined above 
is a pointwise condition on $\bs\theta$. In fact, any $0$-adapted coframing associated to $(M,\I)$
 is a local section of a $G_0$-structure $\mathcal{G}_0$ on $M$, where $G_0\subset \GL(5,\R)$ is the subgroup generated by matrices of the form
				\[
					g = \left(\begin{array}{ccc}
								a &{\bf 0}&{\bf 0}\\
								{\bf b}_1&A&0\\
								{\bf b}_2&0&B
							\end{array}
						\right), \qquad a\ne 0;~A,B\in \GL(2,\R); ~{\bf b}_1, {\bf b}_2\in \R^2,
				\]
		and
				\[
					 J = \left(\begin{array}{ccc}
						1 &{\bf 0}&{\bf 0}\\
								{\bf 0}&0&I_2\\
								{\bf 0}&I_2&0
						\end{array}
						\right).
				\]	
Two hyperbolic Monge-Amp\`ere systems are equivalent if and only if their corresponding $G_0$-structures are equivalent.\footnote{Two ${G}$-structures $\mathcal{G}$ and $\hat{\mathcal{G}}$ on a manifold $M$ 
are said to be \emph{equivalent} if there exists a diffeomorphism $\phi: \mathcal{G}\rightarrow \hat{\mathcal{G}}$ such that $\phi^*\hat{\bs\omega} = \bs\omega$, where $\bs\omega, \hat{\bs\omega}$ are the tautological 
$1$-forms on $\mathcal{G}, \hat{\mathcal{G}}$, respectively.}

\subsection{Integrable Extensions and B\"acklund Transformations  \\{\rm(c.f. \cite{AF12}, \cite{AF15})}}

    {\Def Let $(M,\I)$ be an exterior differential system. A \emph{rank-$k$ integrable extension} of $(M,\I)$ is an exterior differential system $(N,\J)$ with a submersion $\pi: N\rightarrow M$ that satisfies the condition: for each $p\in N$, there exists an
   open neighborhood $U\subset N$ $(p\in U)$ such that 
   	\begin{enumerate}[(1)]
	\item{on $U$, the differential ideal $\J$ is algebraically generated by the elements of $\pi^*{\I}$  
	together with 1-forms $\theta_1,\ldots,\theta_k\in \Omega^1(U)$, where $k = \dim N - \dim M$;}
	\item{for any $p\in U$, let $F_p$ denote the fiber $\pi^{-1}(\pi(p))$; the $1$-forms $\theta_1,\ldots,\theta_k$ restrict to $T_p F_p$ to be linearly independent.}
	\end{enumerate}
  \label{DefIntExt} }\vskip 2mm		

{\remark In Definition \ref{DefIntExt}, one can understand $\J$ as defining
  		a connection on the bundle $\pi:N\rightarrow M$ that is flat over the integral manifolds of $\I$.
More specifically, Condition $(1)$ implies that, if $S\subset M$ is an integral manifold of $(M,\I)$, then $\J$ restricts to $\pi^{-1}(S)$ to be Frobenius; hence, locally, $\pi^{-1}(S)$ is foliated by integral manifolds of $(N,\J)$. Condition $(2)$ implies that, restricting to any integral manifold of $(N,\J)$, $\pi$
   is an immersion, whose image is an integral manifold of $(M,\I)$. }

   {\Def A \emph{B\"acklund transformation} relating two exterior differential systems, $(M_1,\I_1)$ and $(M_2,\I_2)$, 
   is a quadruple
    $(N,\B;\pi_1,\pi_2)$ where, for each $i\in \{1,2\}$, $\pi_i: N\rightarrow M_i$  makes $(N,\B)$
    an integrable extension of $(M_i,\I_i)$.  Such a B\"acklund transformation is represented by the diagram:
      				\begin{figure*}[h!]
				\begin{center}
					\begin{tikzcd}[column sep=tiny]
							    &   (N,\mathcal{B}) \arrow[ld,"\displaystyle\pi_1",swap]\arrow[rd,"\displaystyle{{\pi}_2}"]  & 		\\
							(M_1,\I_1)  &  &(M_2,{\I}_2)
 					\end{tikzcd}
					\label{Backlund}
				\end{center}	
				\end{figure*}
\label{defBacklund}}\vskip 2mm
		
{\Def In Definition \ref{defBacklund}, if $M_1,M_2$ have the {same dimension}, 
		which is not required in general, then we define
		the \emph{rank} of $(N,\B;\pi_1,\pi_2)$ to be the fiber dimension of either $\pi_1$ or $\pi_2$. If $(M_i,\I_i)$ $(i = 1,2)$ are equivalent exterior differential systems, then $(N,\B;\pi_1,\pi_2)$ is called
		an \emph{auto-B\"acklund transformation} of either $(M_i, \I_i)$. \label{Backlundrank}}\vskip 2mm

	{\Example Let $(N,\B;\pi,\bar \pi)$ be a rank-$1$ B\"acklund transformation relating two hyperbolic Monge-Amp\`ere systems $(M, \I)$ and $(\bar M, \bar \I)$. On some open subsets $U\subset M$
	and $\bar U\subset \bar M$, we can choose $0$-adapted coframings such that
								\[
									\I = \<\eta^0,\eta^1\W\eta^2, \eta^3\W\eta^4\>, \qquad 
									\bar\I=\<\bar\eta^0, \bar\eta^1\W\bar\eta^2,\bar\eta^3\W\bar\eta^4\>.
								\]
	Let $V = \pi^{-1}U\cap \bar\pi^{-1}\bar U$, assumed to be nonempty. 
	It is easy to see that the Cauchy characteristics of the system $\<\pi^*\eta^0\>$ $($defined on $V$$)$ are precisely the fibers of $\pi|_V$; similarly for $\bar\pi|_V$. 
			Thus, it is natural to regard $(N,\B;\pi,\bar \pi)$ as nontrivial if $\pi^*\eta^0$ and $\bar\pi^*\bar\eta^0$ are linearly independent $1$-forms on $N$. 
			In particular, it follows that, on $V$, the differential ideal $\B$ is algebraically generated by $\pi^*\I$ and $\bar\pi^*\bar\eta^0$
			as well as by $\bar\pi^*\bar\I$ and $\pi^*\eta^0$.	\label{MAexample}
	}\vskip 2mm

 {\Def  Given a fiber bundle $\pi:E\rightarrow B$, for any $p\in E$, the \emph{vertical tangent space} of $E$ at $p$ is by definition the kernel of $\pi_*:T_pE\rightarrow T_{\pi(p)}B$.}
  
  {\Def A B\"acklund transformation $(N,\B;\pi_1,\pi_2)$ is said to be \emph{nontrivial} if the two fibrations $\pi_1,\pi_2$ have distinct vertical tangent spaces at each point $p\in N$.} \vskip 2mm

\section{Monge-Amp\`ere Systems and Their First Invariants \label{InvMA}}

Let $(M,\I)$ be a hyperbolic Monge-Amp\`ere system. 
Let $\mathcal{G}_0$ denote the $G_0$-structure on $(M,\I)$ (see Definition \ref{MAadmissible}).
One can reduce (see \cite{BGG}) $\mathcal{G}_0$ to a $G_1$-structure $\mathcal{G}_1$ on which 
the tautological $1$-forms $\omega^0,\omega^1,\ldots,\omega^4$ satisfy the following 
structure equations:
			\begin{align}
				\ed\left(
					\begin{array}{c}
						\omega^0\\
						\omega^1\\
						\omega^2\\
						\omega^3\\
						\omega^4
					\end{array}
				\right) &= -\left(
							\begin{array}{ccccc}
								\phi_0 &0&0&0&0\\
								0&\phi_1&\phi_2&0&0\\
								0&\phi_3&\phi_4&0&0\\
								0&0&0&\phi_5&\phi_6\\
								0&0&0&\phi_7&\phi_8
							\end{array}	
						\right)\W
				\left(
					\begin{array}{c}
						\omega^0\\
						\omega^1\\
						\omega^2\\
						\omega^3\\
						\omega^4
					\end{array}
				\right)  \label{StrEqnMA}\\
				&\qquad\qquad
				+ 	\left(	
				\begin{array}{c}
				\omega^1\W\omega^2+\omega^3\W\omega^4\\
				(V_1+V_5)\omega^0\W\omega^3+(V_2+V_6)\omega^0\W\omega^4\\
				 (V_3+V_7)\omega^0\W\omega^3+(V_4+V_8)\omega^0\W\omega^4\\
				 (V_8-V_4)\omega^0\W\omega^1+(V_2-V_6)\omega^0\W \omega^2\\
				 (V_3-V_7)\omega^0\W\omega^1+(V_5-V_1)\omega^0\W\omega^2
				\end{array}
				\right),\nonumber
			\end{align}
		where $\phi_0  = \phi_1+\phi_4 = \phi_5 + \phi_8$, and $G_1\subset G_0$ is the subgroup generated
		by
		\begin{equation}
			g = \left(\begin{array}{ccc}a&{\bf 0}&{\bf 0}\\
							{\bf 0}&A&0\\
							{\bf 0}&0&B
					\end{array}\right), ~A,B\in \GL(2,\R), ~a = \det(A) = \det(B),
		\end{equation}
		and
		\begin{equation}			
		 J = \left(\begin{array}{ccc}1&{\bf 0}&{\bf 0}\\
							{\bf 0}&0&I_2\\
							{\bf 0}&I_2&0
				\end{array}\right)\in \GL(5,\R).	\label{defmatrixJ}
		\end{equation}  

{\Def Let $(M,\I)$ be a hyperbolic Monge-Amp\`ere system.  A
	 \emph{$1$-adapted coframing}\footnote{This is not to be
	confused with a $1$-adapted coframing in the sense of Definition \ref{rank1oneadapted}.} 
	 of  $(M,\I)$ with domain $U\subset M$ 
	is a section $\bs\eta: U\rightarrow \mathcal{G}_1$.}\vskip 2mm	

Following \cite{BGG}, we introduce the notation\footnote{These $S_i$ are those defined
in \cite{BGG} with the same notation scaled by $1/2$.} 
\begin{equation}
	S_1:=\left(\begin{array}{cc} V_1&V_2\\V_3&V_4\end{array}\right), \quad
	S_2:=\left(\begin{array}{cc} V_5&V_6\\V_7&V_8\end{array}\right).		\label{S1S2def}
\end{equation}
It is shown in \cite{BGG} that
{\Prop Along each fiber of $\mathcal{G}_1$, 
			\begin{equation}
				S_i(u\cdot g) = a A^{-1} S_i(u) B, \quad (i = 1,2)	\label{S1S2trans}
			\end{equation}
for any $g = \diag(a;A;B)$ in the identity component of $G_1$. Moreover,
			\begin{equation}
				S_1(u\cdot J) = \left(\begin{array}{cc} -V_4&V_2\\V_3&-V_1\end{array}\right), \quad S_2(u\cdot J) = \left(\begin{array}{cc} V_8&-V_6\\-V_7&V_5\end{array}\right).	\label{S1S2trans2}
			\end{equation}  {\label{VTrans}}}	

Proposition \ref{VTrans} has a simple interpretation: the matrices $S_1$ and $S_2$
		correspond to two invariant tensors under the $G_1$-action. In fact,
		one can verify that the quadratic form
		\begin{equation}
			\Sigma_1:= V_3~\omega^1\omega^3-V_1~\omega^1\omega^4+V_4~\omega^2\omega^3 - 
						V_2~\omega^2\omega^4
		\end{equation}
		and the $2$-form
		\begin{equation}
			\Sigma_2:=V_7~\omega^1\W\omega^3 - V_5~\omega^2\W\omega^3+V_8~\omega^1\W\omega^4
						 - V_6~\omega^2\W\omega^4
		\end{equation}
		are $G_1$-invariant, which implies that $\Sigma_1, \Sigma_2$ 
		are locally well-defined on $(M,\I)$.

An infinitesimal version of Proposition \ref{VTrans} will be useful: for $i = 1,2$,
			\begin{equation}
				\ed S_i \equiv \left(\begin{array}{cc} \phi_4 &-\phi_2\\-\phi_3 &\phi_1\end{array}\right) S_i +
							 S_i \left(\begin{array}{cc} \phi_5&\phi_6\\\phi_7&\phi_8\end{array}\right) \mod \omega^0,\omega^1,\ldots,\omega^4.
				\label{VonFiber}			 
			\end{equation}

An important class of Monge-Amp\`ere systems are the \emph{Euler-Lagrange systems}. In the classical calculus of variations, an Euler-Lagrange system is a PDE system whose solutions
correspond to the stationary points of a given first-order functional. In \cite{BGG}, it is shown:

{\Prop {\rm (\cite{BGG})} A hyperbolic Monge-Amp\`ere system is locally equivalent to an Euler-Lagrange system if and only if $S_2$ vanishes. \label{EL}}

{\remark Proposition \ref{EL} says that the property of 
being \emph{Euler-Lagrange} is intrinsically defined, that is, it does not depend on the choice of local coordinates. }

\section{$G$-structure Equations for B\"acklund Transformations	\label{rank1Gstr}}

In Definition \ref{defBacklund}, it may appear that $(N,\B;\pi_1,\pi_2)$ being 
a B\"acklund transformation imposes conditions
on all components in this quadruple. However, when it is a nontrivial
rank-$1$ B\"acklund transformation
relating two hyperbolic Monge-Amp\`ere systems, one only needs to 
impose conditions
on the exterior differential system\footnote{To be more precise, $(N,\B)$ is a 
hyperbolic exterior differential system of type $s = 2$ in the sense of 
\cite{BGH1}.} $(N,\B)$, as the following proposition shows.

	{\Prop {\rm \cite{C01}} An exterior differential system $(N^6,\B)$ is a nontrivial rank-$1$ B\"acklund transformation relating two hyperbolic Monge-Amp\`ere systems if and only if, for each $p\in N$, there exists 
	an open neighborhood
	$V\subset N$ $(p\in V)$, a coframing
	 $(\theta^0,\bar\theta^0,\theta^1,\ldots,\theta^4)$ and nonvanishing
	 functions $A_1,\ldots,A_4$ $(A_1A_4\ne A_2A_3)$
		defined on $V$, satisfying the conditions:
		\begin{enumerate}[\quad$(1)$]
		\item{the differential ideal $\B = \<\theta^0, \bar\theta^0, \theta^1\W\theta^2,\theta^3\W\theta^4\>_{\rm alg}$;}
		
		\item{the vector bundles $E_0 = \lbb\theta^0\rbb$,
		 $E_1 = \lbb\theta^0, \theta^1,\theta^2\rbb$ and $E_2 = \lbb\theta^0,\theta^3,\theta^4\rbb$ 
		 are invariant along the flow of $X$  (see Definition \ref{InvariantBundle}), 
		 where $X$ is a nonvanishing vector field on $V$ that 
		 annihilates $\theta^0,\theta^1,\ldots,\theta^4$;}
		\item[\quad$\bar{({2})}$]{the vector bundles $\bar E_0 =\lbb\bar\theta^0\rbb$, $\bar E_1 = \lbb\bar \theta^0, \theta^1,\theta^2\rbb$ and $\bar E_2 = \lbb\bar \theta^0,\theta^3,\theta^4\rbb$
		 are invariant along the flow of $\bar X$, where $\bar X$ is a nonvanishing vector field
		 on $V$ that annihilates $\bar \theta^0,\theta^1,\ldots,\theta^4$;}
			
		\item[\quad$(3)$]{the following congruences hold:
					\begin{align*}
							\ed\theta^0 &\equiv A_1\theta^1\W\theta^2+A_2\theta^3\W\theta^4 \mod \theta^0,\\
							\ed\bar\theta^0& \equiv A_3 \theta^1\W\theta^2 + A_4\theta^3\W\theta^4 \mod\bar\theta^0.
					\end{align*}	}
		\end{enumerate}			
		\label{backlundByItself}
	}

This proposition has the following corollary.
	
 {\corollary Let $(N^6,\B; \pi_1,\pi_2)$ be a nontrivial rank-$1$ B\"acklund transformation relating two 
	 hyperbolic Monge-Amp\`ere systems. A
	  coframing defined on an open subset $V\subset N$ 
	 that satisfies Conditions $(1)$-$(3)$ in Proposition \ref{backlundByItself} can always be arranged
	  to satisfy the extra condition: $A_2 = A_3 = 1$. \label{0adaptedmotivation}}
	\vskip 2mm
	{\it Proof.} This is obtained by scaling $\theta^0$ and $\bar\theta^0$.\qed
	
	{\Def A coframing as concluded in Corollary \ref{0adaptedmotivation} is said to be \emph{$0$-adapted}
		to the B\"acklund transformation $(N,\B)$.}\vskip 2mm

			Given a nontrivial rank-$1$ B\"acklund transformation $(N,\B; \pi_1,\pi_2)$ relating two hyperbolic
			Monge-Amp\`ere systems, one can ask whether its $0$-adapted coframings are precisely the 
			local sections of a $G$-structure on $N$. However, this is not true.
			For example, consider a $0$-adapted coframing 
			$(\theta^0,\bar\theta^0, \theta^1,\ldots,\theta^4)$ 
			defined on an open subset $U\subset N$ with corresponding functions $A_1,A_4$. 
			Let $T: U\rightarrow \GL(6,\R)$ 
			be the transformation: 
								 \begin{equation}
								 	T(p): (\theta^0,\bar\theta^0, \theta^1,\theta^2,\theta^3,\theta^4)\mapsto \(\frac{1}{A_1}\theta^0,\frac{1}{A_4}\bar\theta^0, \theta^3,\theta^4,\theta^1,\theta^2\), \quad\forall p\in U.\label{switchingTransformation}
								 \end{equation}
It is easy to see that the coframing on the right-hand-side is $0$-adapted. However, the same transformation, when applied to a $0$-adapted coframing with corresponding functions $ A_1', A_4'$ that are different from $A_1, A_4$, 
may not result in a $0$-adapted coframing.

			One simple strategy, as taken by \cite{C01}, to avoid this issue is by, in addition to understanding the subbundles $\lbb\theta^0\rbb$ and $\lbb\bar\theta^0\rbb$ as an ordered pair, fixing an order for the pair of subbundles 
			$\lbb\theta^0,\bar\theta^0,\theta^1,\theta^2\rbb$ and $\lbb\theta^0,\bar\theta^0,\theta^3,\theta^4\rbb$. Once this is considered, all local $0$-adapted
			coframings respecting such an ordering are precisely the local sections of a $G$-structure, where $G\subset\GL(6,\R)$ is the Lie subgroup consisting of matrices of the form
						\begin{align}
							g&=\left(\begin{array}{cccc}
									\det( \Be)&0&0&0\\
									0&\det(\A)&0&0\\
									0&0&\A&0\\
									0&0&0&\Be
								\end{array}\right),  \label{formofg}\\\nonumber\\
							& \A = (a_{ij}), ~ \Be = (b_{ij})~\in \GL(2,\R).   \nonumber
						\end{align}

Now let $\mathcal{G}$ denote this $G$-structure on $N$. 
Let $\bs\omega=(\omega^1,\omega^2,\ldots,\omega^6)$ be the tautological $1$-form on $\mathcal{G}$. 
Let $\g$ be the Lie algebra of $G$. 
Using the conditions in Proposition \ref{backlundByItself} and the reproducing property of $\bs\omega$, 
one can show that $\bs \omega$ satisfies the following structure equations, 
recorded from \cite{C01} with a slight change of notation:	
			\begin{align}
				\ed\left(
					\begin{array}{c}
						\omega^1\\
						\omega^2\\
						\omega^3\\
						\omega^4\\
						\omega^5\\
						\omega^6
					\end{array}
				\right)=& -\left(\begin{array}{cccccc}
								\beta_0&0&0&0&0&0\\
								0&\alpha_0&0&0&0&0\\
								0&0&\alpha_1&\alpha_2&0&0\\
								0&0&\alpha_3&\alpha_0-\alpha_1&0&0\\
								0&0&0&0&\beta_1&\beta_2\\
								0&0&0&0&\beta_3&\beta_0-\beta_1
							\end{array}
							\right)\W \left(
					\begin{array}{c}
						\omega^1\\
						\omega^2\\
						\omega^3\\
						\omega^4\\
						\omega^5\\
						\omega^6
					\end{array}
				\right)	\label{rank1streq}\\
						&  + 	
				\left(
					\begin{array}{c}
						A_1(\omega^3-C_1\omega^1)\W(\omega^4-C_2\omega^1)+\omega^5\W\omega^6\\
						\omega^3\W\omega^4+A_4(\omega^5-C_3\omega^2)\W(\omega^6-C_4\omega^2) \\
						B_1\omega^1\W\omega^2+C_1\omega^5\W\omega^6\\
						B_2\omega^1\W\omega^2+C_2\omega^5\W\omega^6\\
						B_3\omega^1\W\omega^2+C_3\omega^3\W\omega^4\\
						B_4\omega^1\W\omega^2+C_4\omega^3\W\omega^4
					\end{array}
				\right),	\nonumber
			\end{align}
			where the matrix 
			in $\alpha$ and $\beta$ is a $\g$-valued $1$-form, called a \emph{pseudo-connection} 
			of $\mathcal{G}$; the second term on the right-hand-side is called the 
			\emph{intrinsic torsion} of $\mathcal{G}$. 
			
			It is easy to see that the intrinsic torsion above, as a map defined on $\mathcal{G}$,
			takes values in a $10$-dimensional representation of $G$ and is $G$-equivariant.
			It is proved in \cite{C01} that
			this representation decomposes into 6 irreducible components, 
			as shown by the following equations, where $u\in \mathcal{G}$ is represented as a column, 
			$g$ is
			as in \eqref{formofg}, and $u\cdot g := g^{-1}u$:
			\begin{equation}\label{rank1Torsiontrans}
				\begin{alignedat}{1}
				A_1(u\cdot g) = \frac{\det(\A)}{\det(\Be)}&A_1(u), \quad A_4(u\cdot g)  = \frac{\det(\Be)}{\det(\A)}A_4(u),\\	
				\left(\begin{array}{c}B_1\\B_2\end{array}\right)(u\cdot g) &= \det(\A\Be)\A^{-1}\left(\begin{array}{c}B_1\\B_2\end{array}\right)(u), \\
					\left(\begin{array}{c}B_3\\B_4\end{array}\right)(u\cdot g) &= \det(\A\Be)\Be^{-1}\left(\begin{array}{c}B_3\\B_4\end{array}\right)(u),\\
				\left(\begin{array}{c}C_1\\C_2\end{array}\right)(u\cdot g)&= \phantom{A}\det(\Be)\A^{-1}\left(\begin{array}{c}C_1\\C_2\end{array}\right)(u), \\
					\left(\begin{array}{c}C_3\\C_4\end{array}\right)(u\cdot g)& = \phantom{A}\det(\A)\Be^{-1}\left(\begin{array}{c}C_3\\C_4\end{array}\right)(u).
				\end{alignedat}
			\end{equation}	
			
		{\Def Let $G$ and $\mathcal{G}$ be as above. 
		The B\"acklund transformation\footnote{To be precise, this is a
		 B\"acklund transformation with an ordered pair of characteristic systems.}
		corresponding to $\mathcal{G}$ is
		said to be \emph{generic} if, at each point $u\in \mathcal{G}$, the intrinsic torsion
		takes values in a $G$-orbit
		with the largest possible dimension.}

\section{An Estimate of Generality \label{generality}}
			
			In this section, we address the problem of generality for {\it generic} rank-$1$ B\"acklund transformations
			relating two hyperbolic Monge-Amp\`ere systems. 
			The main ingredients are a {\it $G$-structure reduction procedure} 
			described in \cite{GardnerEquiv} and a {\it theorem of Cartan} described in \cite{Br14}.
			
			{\lemma Let $(N,\B)$ be a 
			nontrivial rank-$1$ B\"acklund transformation relating two hyperbolic Monge-Amp\`ere
			systems. Let $\mathcal{G}$ be the associated $G$-structure. If $(N,\B)$ is generic, then, 
			at each point $u\in\mathcal{G}$,
			the intrinsic torsion
			takes values in an $8$-dimensional $G$-orbit.	\label{rank1genericlemma}}
			\vskip 2mm
			{\it Proof}. Let 
			\[W_1 := {\rm span}((B_1,B_2), (C_1,C_2)), \quad
			W_2 := {\rm span}((B_3,B_4), (C_3,C_4))\]
			 at each point $u\in \mathcal{G}$.
			By \eqref{rank1Torsiontrans}, the function $A_1A_4$ and
			the dimensions of $W_1$ and $W_2$  
			are all invariant under the $G$-action. Let $T$ denote the intrinsic torsion of $\mathcal{G}$.
			We claim that, for each $u\in \mathcal{G}$, the $G$-orbit of $T(u)$
			is at most $8$-dimensional and that this occurs precisely when $W_1$
			and $W_2$ are both $2$-dimensional. To see why this is true, first note that
			if one of $W_i$ $(i = 1,2)$ has dimension less than $2$ at $u\in \mathcal{G}$,
			then the dimension of $u\cdot G$ is at most $7$-dimensional. 
			If both $W_1,W_2$ have dimension $2$ at $u\in \mathcal{G}$, then it is easy to show
			that there exists a unique $g\in G$ such that, at $u' = u\cdot g$,
					\begin{equation}
						\left(\begin{array}{cc}B_1 &C_1\\B_2&C_2\end{array}\right)
						 = \left(\begin{array}{cc}\epsilon_1 &0\\0&1\end{array}\right),\qquad
						 \left(\begin{array}{cc}B_3&C_3\\B_4&C_4\end{array}\right)
						 = \left(\begin{array}{cc} \epsilon_2 &0\\0&1\end{array}\right),	
						 \label{rank1genericBC}
					\end{equation}
			where $\epsilon_i = \pm 1$ $(i = 1,2)$. This completes the proof.
			\qed
			
			Let $(N,\B)$ be a generic rank-$1$
			B\"acklund transformation relating two hyperbolic Monge-Amp\`ere systems.
			By Lemma \ref{rank1genericlemma}, each point $p\in N$ has a
			connected open neighborhood $U\subset N$ on which a canonical coframing $(\omega^1,
			\omega^2,\ldots,\omega^6)$
			can be determined. Such a coframing satisfies
			the equation \eqref{rank1streq}, where all differential forms are defined on 
			$U$ instead of $\mathcal{G}$, and 
			the equations \eqref{rank1genericBC}, where the sign of each $\epsilon_i$ is determined.
			This motivates the following definition.
			
		{\Def Let $N$ be a $6$-manifold. 
		A coframing $(\omega^1,\omega^2,\ldots,\omega^6)$ defined
			on an open subset $U\subset N$ is said to be \emph{$1$-adapted} $($to
			a generic rank-$1$ B\"acklund transformation relating two hyperbolic Monge-Amp\`ere systems$)$
			if
			there exist $1$-forms $\alpha_i, \beta_i$ $(i = 0,\ldots,3)$ and
			functions $A_1,A_4, B_i, C_i$ $(i = 1,\ldots,4)$ defined on $U$ such that  the equations 
			 \eqref{rank1streq} and \eqref{rank1genericBC} are satisfied.	\label{rank1oneadapted}}
			\vskip 2mm
		
		Now we prove the following main theorem that estimates the generality of generic 
		rank-1 B\"acklund transformations relating two hyperbolic Monge-Amp\`ere systems.	
							
		{\theorem Let $N$ be a $6$-manifold. For each $p\in N$, 
		a $1$-adapted coframing  (Definition \ref{rank1oneadapted}) 
		defined on a small open neighborhood $U\subset N$ of $p$
		 can be uniquely determined, up to diffeomorphism, 
		by specifying at most 6 functions of 3 variables.\label{generalitytheorem}}\vskip 2mm

		{\it Proof}. Let $U\subset N^6$ be a sufficiently small connected open 
		subset. Suppose that
				$\bs\omega = (\omega^1,\omega^2,\ldots,\omega^6)$ is a $1$-adapted
				coframing on $U$ in the sense of Definition \ref{rank1oneadapted}. 
				It follows that there exist functions $P_{ij}$ $(i = 0,\ldots,7; ~j = 1,\ldots,6)$ 
				defined on $U$ such that
				$\bs\omega$ satisfies \eqref{rank1streq} and \eqref{rank1genericBC} with
							\[
								\alpha_i = P_{ij}\omega^j,\quad
								\beta_i = P_{i+4,j}\omega^j	\qquad (i = 0,\ldots,3; ~j = 1,\ldots,6).
							\]
				
				There is a standard method to determine the generality of such a coframing $\bs\omega$ up
				to diffeomorphism (see \cite{Br14}). Our application of such
				a method involves mainly three steps.\\
				
				\noindent {\bf Step 1.} By applying $\ed^2 = 0$ to \eqref{rank1streq},
						we find that $P_{ij}$ are related among themselves and with the coefficients of 
						their exterior derivatives. Repeating this, at a stage, no new relations
						among the $P_{ij}$ arise. 
						
						More explicitly, we can  
				choose $s$ expressions
				$a^\alpha$ $(\alpha = 1,\ldots,s)$ from $P_{ij}$, find $r$ expressions $b^\rho$ $(\rho = 1,\ldots,r)$,
				real analytic functions $F^\alpha_i:\R^{r+s}\rightarrow \R$ and  
				$C^i_{jk}:\R^r\rightarrow \R$ satisfying $C^i_{jk}+C^i_{kj} = 0$,
			   	such that
				\begin{enumerate}[\bf (A)]
					\item{the equation \eqref{rank1streq}, in general, takes the form
						\begin{equation}
							\ed\omega^i 
							= -\frac{1}{2}C^i_{jk}(a)\omega^j\W\omega^k\label{cartanproceq1};
						\end{equation}}
					\item{$\ed a^\alpha$, in general, takes the form:
						\begin{equation}
							\ed a^\alpha  = F^\alpha_i(a,b)\omega^i	\label{cartanproceq2};
						\end{equation}
							moreover, applying $\ed^2 = 0$ to
							\eqref{cartanproceq1} yields identities when we take into account both
							\eqref{cartanproceq1} 
							and \eqref{cartanproceq2};	}
					\item{there exist functions $G^\rho_j:\R^{r+s}\rightarrow \R$ such that applying
							$\ed^2  = 0$ to \eqref{cartanproceq2} 
							yields identities when we replace
							$\ed b^\rho$ by $G^\rho_j \omega^j$ and take into account
							\eqref{cartanproceq1} and \eqref{cartanproceq2}.}
				\end{enumerate}	

				\noindent{\bf Step 2.} 
				For the \emph{tableau of free derivatives} associated to $(F^\alpha_i)$, which
				is a subspace of $\Hom(\R^6,\R^s)$ defined at 
				each point of $\R^{r+s}$, compute its \emph{Cartan characters} 
				(an array of $6$ integers $(s_1,s_2,\ldots,s_6)$) 
				and the dimension $\delta$ of its \emph{first prolongation}. 
				For details, 
				see \cite{Br14}. Moreover, in our case, we verify that $\sum_{i = 1}^6 s_i = r$.
				\\
				
				\noindent{\bf Step 3.} Restricting to a domain $V\subset \R^{r+s}$ 
				where the Cartan characters are constants, 
				compare $s:=\sum_{j = 1}^6 js_j$ with $\delta$. By \emph{Cartan's inequality},
				there are two possibilities: either $s = \delta$ (called the \emph{involutive case}) 
				or $s>\delta$.
				
				In the involutive case,
				one can conclude that
				 (see Theorem $3$ in \cite{Br14}):
				\begin{quote}
					\emph{For any $(a_0,b_0)\in \R^{s+r}$ there exists a coframing $\bs\omega$ and functions 
					 $a = (a^\alpha),b = (b^\rho)$ defined 
					 on an open neighborhood of $\bs0\in \R^6$ that satisfy
					 \eqref{cartanproceq1}, \eqref{cartanproceq2} and $(a(\bs0),b(\bs0)) = (a_0,b_0)$. 
					 Moreover,
					locally, such a coframing can be uniquely determined up to diffeomorphism
					 by specifying $s_k$ functions of
					$k$ variables, where $s_k$ is the last nonzero Cartan character.}
				\end{quote}
				
				In the non-involutive case, which is the case we encounter, a natural step to take is
				to \emph{prolong} (see \cite{Br14}) the system
				by introducing the derivatives of $b^\rho$, carry out similar steps as the above, and obtain
				new tableaux of free derivatives
				with Cartan characters $(\sigma_1,\sigma_2,\ldots,\sigma_6)$. 
							
				In practice, however, we do not actually prolong, 
				for it is easy to show that, if $s_k$  is the last nonzero character in $(s_1,\ldots,s_6)$, then 
				$\sigma_j = 0$ $(j>k)$ and $\sigma_k\le s_k$. 
				Using this and the \emph{Cartan-Kuranishi Theorem} (\cite{BCG}),
				one can already conclude that the `generality' of $1$-adapted coframings
				is {\it bounded from above} by $s_k$ functions of $k$ variables.
				(In our case, $k = 3$ and $s_k = 6$.)
				
				For the details of carrying out the steps above, see
				Appendix \ref{App:genThm}. 
				Most of our calculations are performed 
				using Maple\texttrademark. \qed

{\remark		 Clearly, two $1$-adapted coframings that are equivalent under
			a diffeomorphism
			correspond to equivalent B\"acklund transformations. Because of this,
			the upper bound
			for the `generality' of $1$-adapted coframings in Theorem \ref{generalitytheorem} 
			applies to the `generality' of 
			generic rank-$1$ B\"acklund transformations relating two hyperbolic Monge-Amp\`ere systems.}

{\corollary There exist hyperbolic Monge-Amp\`ere systems that are not related to any hyperbolic Monge-Amp\`ere system by a generic rank-$1$ B\"acklund transformation. \label{generalityCor}}
\vskip 2mm
{\it Proof}.  A hyperbolic Monge-Amp\`ere system, up to contact equivalence, can be uniquely 
		determined by specifying $3$ functions of $5$ variables.\footnote{
		One can also apply the same method used in the proof of Theorem \ref{generalitytheorem}
		to verify the stronger statement: Locally, a hyperbolic Euler-Lagrange system, which is Monge-Amp\`ere,
		can be determined uniquely by specifying $1$ function of $5$ variables. (This is not surprising, 
		as,
		in our case,
		a Lagrangian is a function depending on $5$ variables.)}
		
		Note that a generic rank-$1$ B\"acklund transformation in consideration 
		completely determines the two underlying hyperbolic Monge-Amp\`ere systems (up to equivalence).
		The conclusion follows. \qed

\section{Classifications and Examples in Higher Cohomogeneity\label{rk1example}}

			Following the discussion in the previous section, let $U\subset N^6$ be a sufficiently small
			connected
			open subset. Let $\bs\omega$ be a $1$-adapted coframing defined on $U$ 
			in the sense of Definition \ref{rank1oneadapted}. One can ask, \emph{when we specify several
			structure invariants, can we classify the corresponding B\"acklund transformations, if any?}
						
			In the rest of this section, we consider the case when 
			$\epsilon_1 = \epsilon_2 = 1$ in \eqref{rank1genericBC}, and when
			$A_1$ and $A_4$ (or $P_{81}$ and $P_{84}$ in the new notation) in \eqref{rank1streq}
			are specified to be $A_1 = 1$ and $A_4 = -1$.
			
			The following procedure is similar to that in Appendix \ref{App:genThm}. 
			All calculations below are performed
			using Maple\texttrademark.
			
			First, all coefficients in \eqref{rank1streq} are expressed in terms of the remaining
			$40$ $P_{ij}$. Defining their derivatives $P_{ijk}$ by 
				\[
					\ed P_{ij} = P_{ijk}\omega^k
				\]
			and applying the identity $\ed^2 = 0$ to
			 \eqref{rank1streq}, we obtain a system of $106$ polynomial equations in $P_{ij}$ and $P_{ijk}$,
			 which implies that
			\begin{equation*}
				\begin{alignedat}{5}
				&P_{01} = P_{41}, &&{\quad} P_{02} = P_{42}, &&{\quad}  P_{04} = P_{44}, &&{\quad}  P_{06} = P_{46}, &&{\quad}  P_{11} = 0,\\
				& P_{12} = 0,	&&{\quad} P_{16} = 2P_{46}, &&{\quad}  P_{21} = -1, &&{\quad}  P_{22} = -1, &&{\quad}  P_{23} = 0, \\
				& P_{35} = 0, 	&&{\quad}  P_{36} = -1, &&{\quad}  P_{51} = 0, &&{\quad}  P_{52}=0, &&{\quad}  P_{54} = 2P_{44}, \\
				& P_{61} = 1, 	&&{\quad}  P_{62}= -1, &&{\quad}  P_{65} = 0, &&{\quad}  P_{73} = 0, &&{\quad}  P_{74} = -1.
				  \end{alignedat}
			\end{equation*} 

			Using these relations and repeating the steps above, we obtain a system of $88$ equations,
			which implies that
			\begin{align*}
				P_{25}& = -2(P_{41}+P_{44}-P_{42}), \\
				 P_{26} &= P_{64}, \\
				  P_{63} &= -2(P_{41}+P_{42}+P_{46}).
			\end{align*}
			
			Using these and repeating, we obtain a system of $86$ equations for the $17$ $P_{ij}$ remaining
			and $80$ of their $102$ derivatives. This system implies that
			\[
			     P_{31} = P_{32}, \quad P_{41} = -P_{44} - P_{46}, \quad P_{42} = P_{44} - P_{46}, \quad P_{71} = -P_{72}.
			\]

			Using these and repeating, we obtain a system of $85$ equations for the $13$ $P_{ij}$ remaining
			and $64$ of their derivatives. This system implies that
			\begin{equation*}
				\begin{alignedat}{5}
				&P_{03} = -1, &&{\quad}  P_{05} =0, &&{\quad}  P_{15} = 0, &&{\quad}  P_{32} = P_{72}, &&{\quad}  P_{34} = -1, \\ 
				&P_{43} = 0, &&{\quad}  P_{45} = 1, &&{\quad}  P_{53} = 0, &&{\quad}  P_{76} = 1.&&	
				\end{alignedat}
			\end{equation*}
			
			Using these and repeating, we obtain a system of $61$ equations for 
			$P_{72}, P_{44}$, $P_{46}, P_{64}$ and $22$ of their derivatives. 
			Solving this system leads to the two cases below.

			\subsection*{\bf Case $\bf 1$: $P_{72}\ne 0$} 
			
			In this case, we have 
					\[
						P_{44} = 0,\quad P_{46} = 0.
					\]
				Using these and applying $\ed^2=0$ to the structure equations, we find that
					\[
						P_{64} = \frac{1}{P_{72}}.
					\]
				Using this and repeating, we find that 
					\[
						\ed(P_{72}) = 0.
					\]	
				It follows that the only primary invariant remaining, $P_{72}$, is a nonzero constant.
				The structure equations read
				\begin{align*}
					\ed\omega^1& = \omega^1\W(\omega^3+\omega^5) +\omega^3\W\omega^4+\omega^5\W\omega^6,\\
					\ed\omega^2& = -\omega^2\W(\omega^3+\omega^5)+\omega^3\W\omega^4-\omega^5\W\omega^6,\\
					\ed\omega^3& = \(\omega^1+\omega^2-\frac{1}{P_{72}}\omega^6\)\W\omega^4+\omega^1\W\omega^2,\\
					\ed\omega^4& = -(P_{72}\omega^1+P_{72}\omega^2 - \omega^6)\W\omega^3+\omega^5\W\omega^6,\\
					\ed\omega^5& = -\(\omega^1-\omega^2+\frac{1}{P_{72}}\omega^4\)\W\omega^6+\omega^1\W\omega^2,\\
					\ed\omega^6& = (P_{72}\omega^1-P_{72}\omega^2+\omega^4)\W\omega^5+\omega^3\W\omega^4.
				\end{align*}
				This, after the transformation 
				\begin{align*}(\omega^1,\omega^2,&\omega^3,\omega^4,\omega^5,\omega^6)\mapsto\\
				&(\sqrt{|P_{72}|}\omega^1,\sqrt{|P_{72}|}\omega^2,\omega^3,\sqrt{|P_{72}|}\omega^4,\omega^5,\sqrt{|P_{72}|}\omega^6),
				\end{align*}		
				can be readily seen to belong to Case 3D 
				in Clelland's classification (see \cite{C01}). According to
				\cite{C01},
				 if $P_{72}<0$, then $(N,\B)$ is a homogeneous B\"acklund transformation relating time-like
				  surfaces of the constant mean curvature  
				 \[
				 	H =-\frac{ P_{72}}{\sqrt{(P_{72})^2+1}}
				\] 
	in $\mathbb{H}^{2,1}$; if
$P_{72}>0$, then $(N,\B)$ is a homogeneous B\"acklund transformation relating certain surfaces in a $5$-dimensional quotient space of the Lie group $\SO^*(4)$.
\\

\subsection*{\bf Case $\bf 2$: $P_{72} = 0$}

In this case, all coefficients in \eqref{rank1streq} are expressed in terms of $P_{44}, P_{46}$ and $P_{64}$.
Applying $\ed^2 = 0$ to the structure equations, no new relations between $P_{44}, P_{46}$ and $P_{64}$ arise. Furthermore, $16$ of the $18$ derivatives of $P_{44}, P_{46}$ and $P_{64}$ are expressed in terms of these three
invariants; two derivatives, $P_{644}$ and $P_{646}$, are free. 

It is easy to check that Theorem 3 in \cite{Br14} applies to $\bs\omega$, the expressions
$a = (P_{44}, P_{46}, P_{64})$, $b = (P_{644}, P_{646})$, and the functions
$C^i_{jk}$ and $F^\alpha_i$ determined during the calculation above. 
The corresponding \emph{tableaux of free
derivatives} is \emph{involutive} with Cartan characters $(1,1,0,0,0,0)$. We have thus proved the following theorem.

\begin{theorem}{
      Locally, a 
      generic rank-$1$ B\"acklund transformation $(N,\B)$ relating two hyperbolic Monge-Amp\`ere systems 
      with its $1$-adapted coframing 
      satisfying $\epsilon_1 = \epsilon_2 = 1$ and  $A_1 = - A_4 = 1$ can be uniquely determined by specifying 
      $1$ function of $2$ variables. }
 \end{theorem} 
      \vskip 2mm
      
We can study {\bf Case 2} in greater detail. For convenience, we introduce the following new notation:
		\begin{equation}
			\begin{array}{c}
			R:=P_{44}+P_{46}, \quad S:=P_{44}-P_{46}, \quad T:=P_{64};\\[1em]
						 T_{4}:=P_{644}, \quad T_6:=P_{646}.
			\end{array}	\label{RST2P}
		\end{equation}
In this new notation, $\bs\omega$ satisfies \eqref{rank1streq} and \eqref{rank1genericBC} where
		\begin{equation}\label{alphabeta2omega}
			\begin{alignedat}{1}
			\alpha_0& = -R\omega^1+S\omega^2-\omega^3+\tfrac{1}{2}(R+S)\omega^4+\tfrac{1}{2}(R-S)\omega^6,\\
			\beta_0& = -R\omega^1+S\omega^2+\tfrac{1}{2}(R+S)\omega^4+\omega^5+\tfrac{1}{2}(R-S)\omega^6,\\
			\alpha_1& = (R-S)\omega^6,\\
			\alpha_2&= (R+S)\omega^5 + T\omega^6 -\omega^1 - \omega^2,\\
			\alpha_3&=  -\omega^4-\omega^6,\\
			\beta_1&=(R+S)\omega^4,\\
			\beta_2&=(R-S)\omega^3+ T\omega^4+\omega^1-\omega^2,\\
			\beta_3&=-\omega^4+\omega^6,
		\end{alignedat}
		\end{equation}
and $\epsilon_1 = \epsilon_2 = 1$, $A_1 = - A_4 = 1$.		
		
Moreover, the exterior derivatives of $R,S$ and $T$ are
		\begin{equation}\label{dRdSdT}
			\begin{alignedat}{1}		
			\ed R& = -R^2\omega^1+(RS-1)\omega^2+R\omega^3 + \tfrac{1}{2}(R^2+RS-1)\omega^4\\
			&+R\omega^5 + \tfrac{1}{2}(R^2-RS+1)\omega^6, \\
			\ed S& = -(RS-1)\omega^1+S^2\omega^2-S\omega^3+\tfrac{1}{2}(S^2+SR-1)\omega^4\\
			&-S\omega^5+\tfrac{1}{2}(SR-S^2-1)\omega^6,\\
			\ed T& = 2(-RT+S)\omega^1+2(ST-R)\omega^2+2(R^2-S^2)(\omega^3+\omega^5)\\
			&+T_4\omega^4+ T_6\omega^6.
		\end{alignedat}
		\end{equation}

Using these equations, we can study the symmetry of the corresponding B\"acklund transformation
$(N,\B)$. Let $U\subset N$ be the domain of a $1$-adapted coframing $\bs\omega$. Let the map
$\Phi: U\rightarrow \R^3$
be defined by
			\[
				\Phi(p) = (R(p),S(p),T(p)).
			\]

{\lemma The map $\Phi$ can never have rank $0$. Moreover, it has
		\begin{itemize}
			\item{rank $1$ if and only if $2RS=1$ and $T = R^2+S^2$;}
			\item{rank $2$ if and only if it does not have rank $1$ and satisfies either 
				\begin{enumerate}[{\rm (1)}]
					\item{$2RS=1$} or 
					\item{$T_4 = (R+S)(T-1)$ and $T_6 = (R-S)(T+1)$;}
				\end{enumerate}	}
			\item{rank $3$ if and only if it does not have rank $1$ or $2$.}
		\end{itemize}	
\label{variousCohomog}}

	{\it Proof.}  In $\ed R\W \ed S\W \ed T$, the coefficients of $\omega^i\W\omega^j\W\omega^k$ are polynomials in $R,S,T,T_4$ and $T_6$. These coefficients have the common factor $2RS-1$. Calculating with Maple\texttrademark, we find that the coefficients of  $\omega^i\W\omega^j\W\omega^k$ in
	$(2RS-1)^{-1}\ed R\W \ed S\W \ed T$ all vanish if and only if $T_4 = (R+S)(T-1)$ and $T_6 = (R-S)(T+1)$. This justifies the conditions for having rank $2$. The condition for `rank-$1$' can be obtained by setting the coefficients of $\omega^i\W\omega^j$ in $\ed R\W \ed S$, $\ed R\W \ed T$ and $\ed S\W \ed T$ to be all zero. By \eqref{dRdSdT}, it is clear that $\ed R\ne 0$ everywhere; hence, $\Phi$ cannot have rank $0$. \qed

{\Def In the current case, if the corresponding B\"acklund transformation $(U,\B)$ has a symmetry whose orbits
are of dimension $6-k$, 
then it is said to have \emph{cohomogeneity $k$}.\\


\noindent{\bf Case of cohomogeneity-$1$.} This occurs precisely when $\rank(\Phi) = 1$. 
In this case, locally $\Phi$ is a submersion to either branch of the curve in $\R^3$ defined by $2RS = 1$ and $T = R^2+S^2$.   
Expressing $S$ and $T$ in terms of $R$, we have, on $U\subset N$,
		\begin{equation}\label{eqndR}
			\begin{alignedat}{1}
			\ed R &= -R^2\omega^1-\frac{1}{2}\omega^2+R(\omega^3+\omega^5)\\
					&+\frac{1}{4}(2R^2-1)\omega^4+\frac{1}{4}(2R^2+1)\omega^6. 
			\end{alignedat}		
		\end{equation}

It is clear that $\ed R$ is nowhere vanishing. Since $R$ is the only invariant, 
each constant value of $R$ determines a 5-dimensional
submanifold $N_R\subset N$, which has a Lie group structure. The Lie group structure can be determined by 
setting the right-hand-side of the equation \eqref{eqndR} to be zero, obtaining, say,
		\[
			\omega^1 = -\frac{1}{2R^2}\omega^2+\frac{1}{R}\omega^3+
							\frac{2R^2-1}{4R^2}\omega^4+\frac{1}{R}\omega^5
							+\frac{2R^2+1}{4R^2}\omega^6,
		\]
then substituting this into the structure equations, 
yielding equations of $\ed\omega^i$ $(i = 2,\ldots,6)$, expressed
in terms of $\omega^2,\ldots,\omega^6$ alone.	 These are the structure equations on each $N_R$. 
Let $X_1,X_2,\ldots,X_5$ be the vector fields tangent to $N_R$ and
dual to $\omega^2,\ldots,\omega^6$, such that $\omega^i(X_{j}) = \delta^i_{j+1}$ $(i -1,j= 1,2,\ldots,5)$.
We obtain the Lie bracket relations:
	\begin{align*}
		[X_1,X_2]& = 2X_1+\frac{1}{R}(X_2+X_4),\\
		 [X_1,X_3]& = -\frac{1}{2R}X_1
								-\frac{(2R^2-1)}{4R^2}(X_2-X_4)+\frac{1}{R}X_3,\\
		[X_1,X_4]& = 2X_1+\frac{1}{R}(X_2+X_4), \\
		 [X_1,X_5]& = \frac{1}{2R}X_1
									+\frac{2R^2+1}{4R^2}(X_2-X_4)+\frac{1}{R}X_5,\\
		[X_2,X_3]& = -X_1-\frac{1}{R}X_2-X_3-X_5, \\[0.6em]
		[X_2,X_4 ] &= 0,\\[0.3em]
		[X_2,X_5]& = -\frac{2R^2-1}{2R}X_2+X_3+\frac{2R^2+1}{2R}X_4-X_5,\\
		[X_3,X_4]& = -\frac{2R^2-1}{2R}X_2+X_3+\frac{2R^2+1}{2R}X_4-X_5,\\
		[X_3,X_5]& = \left(-R^2+\frac{1}{2}\right)X_2+RX_3+\left(R^2+\frac{1}{2}\right)X_4-RX_5,\\
		[X_4,X_5]&  = X_1-X_3+\frac{1}{R}X_4-X_5.	
	\end{align*}

Using these relations, it can be verified that $X_i$ $(i = 1,\ldots,5)$ generate a 5-dimensional 
Lie algebra that is solvable
but not nilpotent. The derived series has dimensions $(5,3,1,0,\ldots)$.
In fact, after introducing the following new basis:
	\begin{align*}
		\e_1 & = RX_1, \\[0.5em]
		 \e_2 &= -\frac{1}{2}(X_2-X_4), \\
		\e_3& = -\frac{1}{2R}X_1-\frac{2R^2-1}{4R^2}(X_2-X_4)+\frac{1}{R}X_3,\\
		\e_4& = \frac{1}{2R}X_1+\frac{2R^2+1}{4R^2}(X_2-X_4)+\frac{1}{R}X_5,\\
		  \e_5 &= \frac{1}{R}X_1+\frac{1}{2R^2}(X_2+X_4),
	\end{align*}
we obtain the Lie bracket relations:
	\begin{equation*}
		\begin{alignedat}{3}
		&[\e_1,\e_3] = \e_3, &&{\quad} [\e_1,\e_4] = \e_4, && {\quad}[\e_1,\e_5] = 2\e_5,\\
		&[\e_2,\e_3] = \e_4, &&{\quad} [\e_2,\e_4] = -\e_3, &&{\quad} [\e_3,\e_4] = \e_5,
		\end{alignedat}
	\end{equation*}	
with all $[\e_i,\e_j]$ that are not on this list being zero. An equivalent way of writing these relations is:
	\begin{equation*}
		\begin{alignedat}{3}
			&[\e_1+i\e_2, &~\e_1-i\e_2] &= 0,\\
			&[\e_1+i\e_2,&~\e_3+i\e_4] &= 2(\e_3+i\e_4), \\
			 &[\e_1-i\e_2, &~\e_3+i\e_4] &= 0, \\
			 &[\e_1+i\e_2,&~\e_5] &= 2\e_5,\\
			 &[\e_3+i\e_4,&~\e_3-i\e_4] &= -2i\e_5, \\
			&[\e_3+i\e_4,&~\e_5] &= 0.
	         \end{alignedat}
	\end{equation*}	
Now, it is easy to see that the Lie algebra
$\bigoplus_{i = 1}^5 \R\e_i$ is isomorphic to the Lie algebra generated by
the real and imaginary parts of the vector fields 
		\[
			\partial_w, \quad e^{2w}(\partial_z + i\bar z \partial_\lambda),  \quad e^{2(w+\bar w)}\partial_\lambda
		\]
on $\R\times \C^2$ with coordinates $(\lambda;z,w)$. In fact, an isomorphism is induced by the correspondence
		\[
			\e_1+i\e_2\mapsto \partial_w, \quad \e_3+i\e_4\mapsto e^{2w}(\partial_z + i\bar z \partial_\lambda), \quad \e_5\mapsto e^{2(w+\bar w)}\partial_\lambda.
		\]
\vskip 5mm
	Next, we describe the hyperbolic Monge-Amp\`ere systems
	related by the B\"acklund transformation being considered.
	
{\Prop  A B\"acklund transformation in the current (cohomogeneity-$1$) case is an auto-B\"acklund
transformation of a homogeneous Euler-Lagrange system.\label{Cohomog1Prop}}\vskip 2mm

{\it Proof}. This proof is in two parts. First, we show that the underlying two hyperbolic Monge-Amp\`ere systems
		are equivalent and are homogeneous. Second, by computing their local invariants, 
		we verify that they are hyperbolic 
		\emph{Euler-Lagrange
		systems} in the sense of \cite{BGG}.

Using the structure equations on $U\subset N$, if we let $(\theta^0,\theta^1,\ldots,\theta^4)$ be either
		\begin{equation}
			\left(S\omega^1, -R(\omega^1-\omega^4)+\omega^3, S\omega^4, 
			-R(\omega^1-\omega^6)+\omega^5, S\omega^6\right)	\label{cohomog1MA1}
		\end{equation}
		or
		\begin{equation}
			\left(-R\omega^2, S(\omega^2+\omega^4)-\omega^3, R\omega^4, S(\omega^2-\omega^6)-\omega^5, -R\omega^6\right),									\label{cohomog1MA2}
		\end{equation}
		then we can verify (with Maple\texttrademark) 
		that  $\theta^i$ $(i = 0,\ldots,4)$, in both cases, satisfy the same structure equations:
		\begin{align}
			\ed\theta^0& = -2(\theta^1+\theta^3)\W\theta^0+\theta^1\W\theta^2+\theta^3\W\theta^4,	\nonumber\\
			\ed\theta^1& = -\theta^1\W\theta^4 +\theta^2\W\theta^4 +\theta^2\W\theta^3,	\nonumber\\
			\ed\theta^2& =  -\theta^1\W\theta^4 +\theta^2\W\theta^4 +\theta^2\W\theta^3-\theta^1\W\theta^2+\theta^3\W\theta^4,	\label{Cohomog1MAstreq}\\
			\ed\theta^3& =  {\phantom{-}}\theta^1\W\theta^4-\theta^2\W\theta^4 -\theta^2\W\theta^3,	\nonumber\\
			\ed\theta^4& =  - \theta^1\W\theta^4+\theta^2\W\theta^4 +\theta^2\W\theta^3+\theta^1\W\theta^2 - \theta^3\W\theta^4.	\nonumber
		\end{align}
		It is easy to verify that \eqref{Cohomog1MAstreq} are the structure equations on a 
		$5$-manifold $M$ with a hyperbolic Monge-Amp\`ere ideal 
		$\I =\<\theta^0,\theta^1\W\theta^2,
		\theta^3\W\theta^4\>$.
		It follows that the expressions \eqref{cohomog1MA1}
		and \eqref{cohomog1MA2} correspond to the pull-back of $\theta^i$ under two distinct 
		submersions $\pi_1,\pi_2: U\rightarrow M$. It is easy to see that $(U,\B;\pi_1,\pi_2)$ is an
		auto-B\"acklund transformation of the system $(M,\I)$.
		Moreover, $(M,\I)$ is \emph{homogeneous}, since all 
		coefficients in \eqref{Cohomog1MAstreq} are constants.
		
		Next, we verify that $(M,\I)$ is a \emph{hyperbolic Euler-Lagrange system}. In \cite{BGG}, it is proved that
		a hyperbolic Monge-Amp\`ere system is Euler-Lagrange 
		if and only if the
		invariant tensor
		$S_2$ vanishes
		(see Proposition \ref{EL}).
		To compute $S_2$ in the current case, we choose a new coframing $\bs\eta = (\eta^i)$ and $1$-forms 
		$(\phi_\alpha)$ below: 
		\begin{align*}
				(\eta^0,\eta^1,\eta^2,\eta^3,\eta^4) = (\sqrt{2}\theta^0, &~\sqrt{2}\theta^1,~\theta^1+\theta^2-\theta^0,\\
					& \theta^3+\theta^4-\theta^0, ~\sqrt{2}(\theta^4-\theta^0)), 
		\end{align*}
		\begin{equation}\label{etasandphis}
			\begin{alignedat}{2}	
				\phi_0& = \frac{1}{\sqrt{2}}\eta^1 + \eta^3-\frac{1}{\sqrt{2}}\eta^4, \\
				\phi_1 &= -\sqrt{2}\eta^0 - \eta^3 - \frac{1}{\sqrt{2}}\eta^4,&&{\quad}
				\phi_2 = \eta^0+\sqrt{2}\eta^3,\\
				 \phi_3 &= -2\eta^0+\eta^1-\frac{1}{\sqrt{2}}\eta^2-\sqrt{2}\eta^3-\eta^4,&&{\quad}
				\phi_4  = \phi_0-\phi_1, \\
				 \phi_5 &= -\eta^3-\frac{1}{\sqrt{2}}\eta^4, &&{\quad}
				\phi_6 = \frac{1}{\sqrt{2}} \eta^3,\\
				\phi_7 &= -\eta^0+\eta^1-\sqrt{2}\eta^2-\sqrt{2}\eta^3,&&{\quad}  \phi_8 = \phi_0-\phi_5.
			\end{alignedat}
		\end{equation}
		
		These $\eta^i$ $(i = 0,1,\ldots,4)$ and $\phi_\alpha$ $(\alpha = 1,\ldots,8)$ are chosen such that they
		satisfy the structure equations \eqref{StrEqnMA}, and such that $S_1$ and $S_2$
		are as simple as possible. One can verify that, under this choice,
				\[
					S_1 = \left(\begin{array}{cc}1&0\\0&1\end{array}\right),  \quad  S_2 = \bs0.
				\]
		This completes the proof. \qed
	
		{\remark\label{nonFGordon}
			The Monge-Amp\`ere system $(M,\I)$ considered in Proposition \ref{Cohomog1Prop} cannot,
			by a contact transformation, be put in the form
			\[
				z_{xy} = f(x,y,z,z_x,z_y).
			\]
			This is because $\rank(S_1) = 2$ and $S_2 = \bs0$, which implies that 
			neither of the two 
			characteristic systems of $(M,\I)$ contains a rank-$1$ integrable subsystem.
		}\vskip 2mm
		
		Now, one may wonder whether the homogeneous Monge-Amp\`ere system $(M,\I)$ 
		considered in
		Proposition \ref{Cohomog1Prop} has a symmetry of dimension greater than $5$. 
		Using the \emph{method of equivalence}, we prove that it is not the case.

		{\Prop The hyperbolic Euler-Lagrange system in
		Proposition \ref{Cohomog1Prop} has a
		symmetry of dimension $5$. In addition, any such symmetry is induced from
		a symmetry of the B\"acklund transformation $(N,\B)$. \label{Cohomog1SymmetryProp}}\vskip 2mm					
	{\it Proof.} Let $(M,\I)$ denote the Euler-Lagrange system being considered. To show that $(M,\I)$ has 
			a $5$-dimensional symmetry, it suffices to show that there is a canonical way to determine
			a local coframing on $M$. This can be achieved by applying the \emph{method of equivalence}.
			For details, see Appendix \ref{App:ELSym}.
			
			By \eqref{cohomog1MA1} and \eqref{cohomog1MA2}, it is easy to see that 
			the fibers of $\pi_i: N\rightarrow M$ $(i = 1,2)$ are everywhere transversal to 
			the level sets of the functions $R$. The second half of the statement follows.\qed
			\\

		To end the discussion of the cohomogeneity-$1$ case, we integrate the structure equations \eqref{Cohomog1MAstreq} to express the corresponding hyperbolic Euler-Lagrange system
		in local coordinates. 
		
{\Prop \label{Cohomog1PDEProp}	The hyperbolic Monge-Amp\`ere system $(M,\I)$ with the differential ideal 
$\I = \<\theta^0,\theta^1\W\theta^2, \theta^3\W\theta^4\>$, where $\theta^i$ satisfy 
\eqref{Cohomog1MAstreq}, is equivalent to the following hyperbolic Monge-Amp\`ere PDE up to a contact transformation:
		\begin{equation}
			(A^2-B^2)(z_{xx} - z_{yy}) + 4ABz_{xy} = 0,	\label{Cohomog1PDE}
		\end{equation}
where $A = 2z_x+y$ and $B = 2z_y - x$.}\vskip 2mm
	
{\it Proof}. Let $U\subset M$ be a domain on which $\theta^i$ are defined.
			One can verify, using the structure equations \eqref{Cohomog1MAstreq}, that
			the $1$-forms $\theta^1-\theta^2-\theta^3-\theta^4$ and $\theta^1+\theta^3$ are closed.
			Hence, by shrinking $U$ if needed, there exist functions $P,Q$ defined on $U$ such that
				\[
					\ed P  = -(\theta^1+\theta^3), \quad \ed Q = \theta^1-\theta^2-\theta^3-\theta^4.
				\] 		
			Moreover, if we let $\Theta = \theta^2+i\theta^4$, then, by a straightforward calculation,
			we obtain
				\[
					\ed \Theta = \ed (P+iQ)\W \Theta.
				\]
			It follows that there exist functions $X,Y$ on $U$ such that
				\[
					\Theta = e^{P+iQ}\ed(X+iY).
				\]
			Equivalently, we have
				\begin{align*}
					\theta^2 &= e^P(\cos Q ~\ed X - \sin Q~  \ed Y), \\
					 \theta^4 &= e^P(\sin Q~ \ed X+\cos Q~ \ed Y).
				\end{align*}
			Now we can express $\theta^1,\ldots,\theta^4$ completely in terms of the functions $X,Y,P,Q$.	
			By verifying the equality
				\begin{equation}
					\ed (e^{-2P}\theta^0) = e^{-2P}(\theta^1\W\theta^2+\theta^3\W\theta^4), \label{derivt0}
				\end{equation}
			we notice that the right-hand-side of \eqref{derivt0} is a symplectic form on a $4$-manifold on which
			 $\theta^1,\ldots,\theta^4$ are well-defined, by \eqref{Cohomog1MAstreq}.
			Thus, by the theorem of Darboux, 
			 locally there exist functions $x,y,p,q$ such that the right-hand-side of \eqref{derivt0} 
			is equal to $\ed x\W \ed p+\ed y\W \ed q$.
			
			 In fact, in the $XYPQ$-coordinates, the right-hand-side of \eqref{derivt0} is equal to
			\begin{align*}
				&\ed\left(\frac{e^{-P}}{2}(\cos Q +\sin Q)+\frac{Y}{2}\right)\W \ed X\\
				 &\qquad\qquad	+ \ed\left(\frac{e^{-P}}{2}(\cos Q-\sin Q)- \frac{X}{2}\right)\W \ed Y.
			\end{align*}
			As a result, we can set 
				\begin{equation*}
					\begin{alignedat}{2}
						&x = X,  &&{\quad} p = -\frac{e^{-P}}{2}(\cos Q +\sin Q)-\frac{Y}{2}, \\
						&y = Y, &&{\quad} q = -\frac{e^{-P}}{2}(\cos Q -\sin Q)+\frac{X}{2},
					\end{alignedat}
				\end{equation*}
			and write 
				 \[
				 	e^{-2P}\theta^0  = \ed z - p\ed x - q\ed y,\]
			 for some function $z$,
			independent of $x,y,p,q$. From these expressions, it is clear that $2p+y$ and $2q - x$ 
			cannot simultaneously vanish. 
			
			Now let $A = 2p+y$ and $B = 2q-x$. 
			We can express $\theta^1\W\theta^2$ in terms of $x,y,z,p,q$:
				\[
					\frac{(A^2-B^2)(\ed p\W \ed y - \ed x \W\ed q) + (A+B)^2 \ed x\W \ed p 					 + (A-B)^2 \ed y\W \ed q}{(A^2+B^2)^2}.
				\]
			Multiplying this expression by $(A^2+B^2)^2$ then subtracting the  result by 
			$(A^2+B^2)(\ed x\W \ed p+\ed y\W \ed q)$, 
			we obtain
				\[
					(A^2-B^2)(\ed p\W \ed y - \ed x\W \ed q)+ 2AB(\ed x\W \ed p + \ed q\W \ed y).
				\]
			The vanishing of this $2$-form on integral surfaces (satisfying the independence condition
			$\ed x\W\ed y \ne 0$)
			implies that $z$ must satisfy the equation $\eqref{Cohomog1PDE}$.
			\qed	
			\\
			

\noindent {\bf Case of cohomogeneity-$2$.} By Lemma \ref{variousCohomog}, 
this case can only occur when $2RS=1$ and $R^2+S^2 = T$ do not both hold and 
either 
\begin{enumerate}[(1)]
	\item{$2RS = 1$} or 
	\item{$T_4 = (R+S)(T-1)$, $T_6 = (R-S)(T+1)$} 
\end{enumerate}
	holds. 
We now focus on the latter case.

{\Prop When $\Phi$ has rank $2$, and when $T_4 = (R+S)(T-1)$ and $T_6 = (R-S)(T+1)$, 
the map $\Phi: N\rightarrow \R^3$ has its image contained in a surface that is defined
by either \[\dfrac{R^2+S^2-T}{2RS-1}\] or its reciprocal being a constant.}
\vskip 2mm

{\it Proof.} First note that $R^2+S^2-T$ and $2RS-1$ cannot be both zero, for this would 
reduce to the cohomogeneity-1 case; hence, the conclusion has meaning. 
To see that this statement is true, note that, in the current case, the pull-back of $\ed R, \ed S$ and $\ed T$ via
$\Phi$ to $N$ are linearly dependent. To be precise, the $1$-form
	\begin{align*}
		\theta =& -2(R^2S-S^3+ST-R)\ed R \\
				&+ 2(R^3-RS^2-RT+S)\ed S+(2RS-1)\ed T
	\end{align*}	
equals to zero when pulled back to $N$. Since the tangent map $\Phi_*$ has rank $2$, this can only occur when
$\theta\W \ed \theta = 0$.
It follows that $\theta$ is integrable. In fact, it is easy to verify that 
the primitives of $(2RS-1)^{-2}\theta$ and $(R^2+S^2-T)^{-2}\theta$
are, respectively, the function \[\dfrac{R^2+S^2-T}{2RS-1}\] and its reciprocal when 
$2RS-1$ and $R^2+S^2-T$ are, respectively, nonzero. This completes the proof.\qed
\\
			
We now study the symmetry of the B\"acklund transformation $(N,\B)$ being considered.			
Let $X_i$ $(i = 1,2,\ldots,6)$ be the vector fields defined on $U\subset N$ that are 
dual to $\omega^i$ $(i = 1,2,\ldots,6)$. 
Using the expressions of $\ed R, \ed S$ and $\ed T$, it is easy to see that the rank-$4$ distribution on $U$ annihilated
by $\ed R,\ed S$ and $\ed T$ is spanned by the vector fields
	\begin{align*}
		Y_1 &= RX_2+SX_1+X_5, & Y_2 = \frac{1}{2}(X_1 - X_2)+X_4,\\
		Y_3 &= X_3 -X_5, & Y_4 = \frac{1}{2}(X_1+X_2)+X_6.
	\end{align*}
These vector fields generate a $4$-dimensional Lie algebra $\mathfrak{l}$ with
	\begin{align*}
		[Y_1,Y_2]& = \frac{R+S}{2}Y_3+Y_4, \\[0.3em]
		 [Y_1,Y_3] &= 0, \\
		  [Y_1,Y_4] &= -Y_2+\frac{R-S}{2}Y_3,\\[0.3em]
		[Y_2,Y_3]& = (R+S)Y_3+2Y_4, \\[-0.1em]
		[Y_2,Y_4] &= Y_1+(R-S)Y_2 +\left(\frac{1}{2}-T\right)Y_3-(R+S)Y_4,\\[0.1em]
		 [Y_3,Y_4] &= 2Y_2 - (R-S)Y_3.\\
	\end{align*}			
It is easy to verify that $2Y_1+Y_3$ belongs to the center of $\mathfrak{l}$. The quotient algebra
$\mathfrak{q} = \mathfrak{l}/\R(2Y_1+Y_3)$, with the basis $\e_1 = [Y_2]$, $\e_2 = [Y_3]$ and $\e_3 =[Y_4]$,
satisfies
	\begin{align*}
		[\e_1, \e_2] &= (R+S)\e_2 + 2\e_3, \\
		 [\e_1,\e_3] &=(R-S)\e_1 -T\e_2-(R+S)\e_3, \\
		  [\e_2,\e_3] &= 2\e_1+(S-R)\e_2.
	\end{align*}
According to the classification of $3$-dimensional Lie algebras, see Lecture 2 in \cite{BryantSymp} for example,
to identify the Lie algebra $\mathfrak{q}$, it suffices to find a normal form of the matrix (note that it is symmetric)
		\[
			C = \left(\begin{array}{ccc}
					2&S-R&0\\
					S-R&T&R+S\\
					0&R+S&2
					\end{array}\right)
		\]	
under the transformation $C\mapsto \det(A^{-1})AC A^{T}$, where $A\in \GL(3,\R)$. Note that
$\det(C) = -2(R^2+S^2-T)$. We have:

{\Prop If $R^2+S^2<T$, then $\mathfrak{q}$ is isomorphic to $\so(3,\R)$. If $R^2+S^2>T$, then 
		$\mathfrak{q}$ is isomorphic to $\sl(2,\R)$. If $R^2+S^2 = T$, then $\mathfrak{q}$ is isomorphic to the 
		solvable Lie algebra with a basis $x_1,x_2,x_3$ satisfying $[x_2,x_3] = x_1$, $[x_3,x_1] = x_2$ and $[x_1,x_2] = 0$}.
\vskip 2mm

	{\it Proof.}  After a transformation of the form above, $C$ can be put in the form
			\[
				C' = \left(\begin{array}{ccc}
					2&0&0\\
					0&T-R^2-S^2&0\\
					0&0&2
					\end{array}\right).
			\]		
			By \cite{BryantSymp}, the conclusion follows immediately. \qed
\\

Now consider the case when $\mathfrak{q}$ is solvable, that is, when $R^2+S^2 = T$.
By the cohomogeneity-$2$ assumption, we must have $2RS\ne 1$. We proceed to identify the 
Monge-Amp\`ere systems related by such a B\"acklund transformation.

If we let  $(\theta^0,\theta^1,\ldots,\theta^4)$ be
			\[
				(S\omega^1, -R(\omega^1-\omega^4)+\omega^3, S\omega^4, -R(\omega^1-\omega^6)+\omega^5, S\omega^6)
			\] 
and let $F$ be defined by
			\[
				F = \frac{2RS-1}{2S^2},
			\]
for which to have meaning we need to restrict to a domain on which $S\ne 0$, then the structure equations on $N$ would imply
		\begin{equation}\label{Cohomog2MAstreq}
			\begin{alignedat}{1}
			\ed\theta^0& = \theta^0\W(2\theta^1-F\theta^2+2\theta^3-F\theta^4)+\theta^1\W\theta^2+\theta^3\W\theta^4,\\
			\ed\theta^1& = -F\theta^0\W(\theta^2+\theta^4) - \theta^1\W\theta^4+\theta^2\W\theta^3+(F+1)\theta^2\W\theta^4,\\
			\ed\theta^2& = -2F\theta^0\W\theta^2-\theta^1\W\theta^2-\theta^1\W\theta^4+\theta^2\W\theta^3\\
			&+(1-F)\theta^2\W\theta^4+\theta^3\W\theta^4,\\
			\ed\theta^3& = F\theta^0\W(\theta^2+\theta^4)+\theta^1\W\theta^4-\theta^2\W\theta^3 +(F-1)\theta^2\W\theta^4, \\
			\ed\theta^4& =-2F\theta^0\W\theta^4+\theta^1\W\theta^2-\theta^1\W\theta^4+\theta^2\W\theta^3\\
			&+(F+1)\theta^2\W\theta^4 - \theta^3\W\theta^4, 
			\end{alignedat}
		\end{equation}
and 
		\begin{equation}
			\ed F = 2F^2(2\theta^0- \theta^2-\theta^4) + 2F(\theta^1+\theta^3), \quad (F\ne 0).\label{Cohomog2F}
		\end{equation}
		
It can be verified that, in the equations \eqref{Cohomog2MAstreq} and \eqref{Cohomog2F},
the exterior derivative of the right-hand-sides are zero, by taking into account these equations themselves.
By the construction of the $\theta^i$ $(i = 0,\ldots,4)$, it follows that \eqref{Cohomog2MAstreq} and \eqref{Cohomog2F}
are the structure equations of one of the Monge-Amp\`ere systems being related by 
the B\"acklund transformation $(N,\B)$. 

On the other hand, it is easy to verify that the transformation
	\[
		(\theta^0,\theta^1,\theta^2,\theta^3,\theta^4; F)\mapsto (-\theta^0, \theta^3,-\theta^4,\theta^1,-\theta^2;-F)
	\]
leaves \eqref{Cohomog2MAstreq} and \eqref{Cohomog2F} unchanged. Thus, by applying
such a transformation, if needed, we can assume that $F>0$.

{\Prop The hyperbolic Monge-Amp\`ere system $(M,\I)$ with the differential ideal 
$\I = \<\theta^0,\theta^1\W\theta^2, \theta^3\W\theta^4\>$, where $\theta^i$ satisfy 
\eqref{Cohomog2MAstreq} and \eqref{Cohomog2F}, corresponds 
to the following hyperbolic Monge-Amp\`ere PDE up to a contact transformation:
		\begin{equation}
			(A^2-B^2)(z_{xx}-z_{yy})+4ABz_{xy} + (A^2+B^2)^2 = 0.	\label{Cohomog2PDE}
		\end{equation}
where $A = z_x-y$, $B = z_y + x$.		\label{Cohomog2PDEProp}}
\vskip 2mm

{\it Proof}. The proof is similar to that of Proposition \ref{Cohomog1PDEProp}. First it is easy to verify that 
	the $1$-forms $F(-2\theta^0+\theta^2+\theta^4) - \theta^1-\theta^3$ and $\theta^1-\theta^2-\theta^3-\theta^4$ are closed. Consequently, locally there exist functions $f, g$ such that 
			\begin{align*}
				\ed f &= F(-2\theta^0+\theta^2+\theta^4) - \theta^1-\theta^3,\\
				 \ed g &= \theta^1-\theta^2-\theta^3-\theta^4.
			\end{align*}
		Now the expression of $\ed F$ can be written as $\ed F = -2F\ed f$. 
		This implies that there exists a constant
		$C>0$ such that $F = Ce^{-2f}$. Using the ambiguity in $f$ (as $f$ is determined up to an additive
		constant), we can arrange that $C = 1$.
		In addition, if we let $\Theta = e^{-f}(\theta^2+i\theta^4)$, 
		it is easy to verify that $\Theta$ is integrable. To be explicit,
			\[
				\ed\Theta = i~ \ed g \W \Theta.
			\]
		Thus, there exist functions $X,Y$ such that  $\Theta = e^{ig}(\ed X+i\ed Y)$. From this we obtain
			\begin{align*}
				\theta^2 &= e^f(\cos g ~\ed X - \sin g ~\ed Y), \\
				 \theta^4 &= e^f(\sin g ~\ed X + \cos g ~\ed Y).
			\end{align*}
		Using these, differentiating $\theta^1+\theta^3$ gives
			\[
				\ed(\theta^1+\theta^3) = 2~\ed X\W \ed Y.
			\]
		This implies that there exists a function $Z$, independent of $X,Y,f,g$, such that 
			\[
				\theta^1+\theta^3 = \ed Z + X\ed Y-Y \ed X.
			\]	
		Now, $\theta^0,\theta^1,\ldots,\theta^4$ 
		can be completely expressed in terms of the functions $X,Y,Z,f,g$.
		In particular,
			\begin{align*}
				-2e^{-2f}\theta^0 =  \ed(Z+f) &- (Y+e^{-f}(\sin g+\cos g)) \ed X \\
										&+(X+e^{-f}(\sin g - \cos g))\ed Y.
			\end{align*}
		If we make the substitution
			\begin{equation*}
				\begin{alignedat}{2}
					&x = X,  &&{\quad} p = e^{-f}(\cos g+\sin g)+Y,\\
					 &y = Y, &&{\quad} q = e^{-f}(\cos g - \sin g)-X,\\[0.4em]
					 &  z = Z+f,&&
				\end{alignedat}
			\end{equation*}	
		the contact form $\theta^0$ is then a nonzero multiple of $\ed z- p\ed x - q \ed y$. The $2$-form
			$\theta^3\W\theta^4$, when each $\ed z$ is replaced by $p\ed x + q\ed y$, can be expressed as
			\begin{equation*}\begin{alignedat}{2}
				\theta^3\W\theta^4\equiv \frac{1}{8}e^{4f} \Big\{&(A^2-B^2)(&&\ed p\W \ed y - \ed x\W \ed q)\\
						+&~(A+B)^2&&\ed q\W \ed y\\
						-&~(A-B)^2 &&\ed x\W \ed p\\
						+&(A^2+B^2)^2&&\ed x\W \ed y\Big\}\qquad 
				\end{alignedat}	\qquad	\mod \theta^0,
			\end{equation*}
		where $A = p-y$, $B = q+x$. Note that, by construction, $A,B$ cannot be simultaneously zero.
		The equation \eqref{Cohomog2PDE} follows.	\qed
		\\
		
In the current case, there remain several obvious questions to investigate. 
{\it  What is the Monge-Amp\`ere system
corresponding to $\<\omega^2,\omega^3\W\omega^4,\omega^5\W\omega^6\>$?
Are the Monge-Amp\`ere systems being B\"acklund-related Euler-Lagrange? Is the B\"acklund transformation
an auto-B\"acklund transformation?} Answers to these questions can be obtained in a similar way as in the
cohomogeneity-$1$ case. We thus have them summarized in the following remark, omitting the 
details of calculation.

{\remark\begin{enumerate}[\bf A.] \item{Whenever $R\ne 0$,
	\[	
		(R\omega^2,S(\omega^2+\omega^4)-\omega^3, -R\omega^4, S(\omega^2-\omega^6)-\omega^5, R\omega^6)
	\]
	form a coframing defined on a $5$-manifold. The system 
	\[\<\omega^2,\omega^3\W\omega^4,\omega^5\W\omega^6\>\] descends to correspond to
	the same equation \eqref{Cohomog2PDE} up to a contact transformation.
	The system $(N,\B)$ is therefore an auto-B\"acklund
	transformation of the equation \eqref{Cohomog2PDE}.}\\
	
	\item{One can verify that the
		 hyperbolic Monge-Amp\`ere system in Proposition \ref{Cohomog2PDEProp}
	is Euler-Lagrange. In fact, by a transformation of the $\theta^i$ $(i = 0,\ldots,4)$, 
	the structure equations \eqref{Cohomog2MAstreq} can be put in the form of \eqref{StrEqnMA}
	with
		\[	S_1 = \left(\begin{array}{cc}1&0\\0&1\end{array}\right),  \quad  S_2 = \bs0.\]
	Applying the procedure described in Appendix \ref{App:ELSym},
	we find the corresponding Monge-Amp\`ere invariants
		\[
			Q_{01}Q_{04}-Q_{02}Q_{03} = -\frac{1}{2}(1+F^2)\ne 0, \quad Q_{00} = 0.
		\]
	Using this and the expression of $\ed F$, 
	one can show $($in a way that is similar to the proof of Proposition \ref{Cohomog1SymmetryProp}$)$ that the underlying Euler-Lagrange system has a symmetry
	of dimension $4$. Such a symmetry is induced from the symmetry of the 
 	B\"acklund transformation $(N,\B)$.}\\
	
	\item{Similar to the reason in Remark \ref{nonFGordon}, \eqref{Cohomog2PDE}
		is not contact equivalent to any hyperbolic Monge-Amp\`ere PDE of the form
			\[
				z_{xy} = f(x,y,z,z_x,z_y).
			\]		}
	\end{enumerate}		
}

\section{Concluding Remarks}

\subsection{Regarding Classification in the Rank-$1$ Case}\

The method outlined in Section \ref{rk1example} can be carried further, for example, by putting weaker restrictions
on the invariants. Which new subclasses of B\"acklund transformations relating two hyperbolic Monge-Amp\`ere
systems can be classified?

Equation \eqref{Cohomog1PDE} admits the trivial solution $z(x,y) = 0$. Which solution do we obtain by applying the 
B\"acklund transformation found in Proposition \ref{Cohomog1Prop} to this trivial solution? Does this give rise to a 
$1$-soliton
solution to \eqref{Cohomog1PDE}?

\subsection{Regarding Monge-Amp\`ere Invariants}\

It is interesting to ask which pairs of hyperbolic Monge-Amp\`ere systems may be related by 
a rank-$1$ B\"acklund transformation. This question can be partially answered by studying how obstructions to the existence of B\"acklund transformations may be expressed in terms of the invariants of the underlying hyperbolic Monge-Amp\`ere systems. Some relevant results are presented in \cite{HuBacklund2}.

\subsection{B\"acklund Transformations of Higher Ranks}\

		Among the examples discussed in \cite{Rogers}, a B\"acklund transformation relating solutions of the \emph{hyperbolic Tzitzeica equation}
is particularly interesting. The \emph{hyperbolic Tzitzeica equation} is the second-order
equation for $h(x,y)$:
				\begin{equation}
					(\ln h)_{xy} = h - h^{-2}.	\label{Tz}
				\end{equation}
This equation was discovered by Tzitzeica in his study of 
\emph{hyperbolic affine spheres} in the affine $3$-space $\mathbb{A}^3$ (see \cite{Tz08} and \cite{Tz09}).		
He found that the system in  $\alpha$, $\beta$ and $h$,
				\begin{equation}
					\left\{\begin{array}{l}
					\alpha_{x} =({h_x\alpha+\lambda\beta}) h^{-1} - \alpha^2,\\[0.8em]
					\alpha_y = \beta_x = h - \alpha\beta, \\	[0.8em]
					\beta_y = ({h_y\beta+\lambda^{-1}\alpha})h^{-1} - \beta^2,	
					\end{array}
					\right.\label{Tztrans}
				\end{equation}
where $\lambda$ is an arbitrary nonzero constant,		
is a B\"acklund transformation
relating solutions of \eqref{Tz}. More explicitly, if $h$ solves the hyperbolic Tzitzeica equation \eqref{Tz},
then, substituting it in the system \eqref{Tztrans}, one obtains a compatible first-order PDE system for $\alpha$ and $\beta$, 
whose solutions 
can be found by solving ODEs; for each solution $(\alpha,\beta)$, the function
					\[
						\bar h = -h + 2\alpha\beta
					\]
also satisfies the hyperbolic Tzitzeica equation \eqref{Tz}. 

Unlike the systems \eqref{CR} and \eqref{SGT}, substituting a solution 
 $h$ of \eqref{Tz} into \eqref{Tztrans} (with fixed $\lambda$) yields a system whose solutions depend on 2 parameters instead of 1. Using our 
 terminology (see Definition \ref{Backlundrank}), 
 one can verify that 
 the system \eqref{Tztrans} corresponds to a \emph{rank}-$2$ B\"acklund transformation.	
		
Furthermore, in \cite{AF15}, it is shown that the following hyperbolic Monge-Amp\`ere equation
	\[
		z_{xy} = \frac{\sqrt{1- z_x^2}\sqrt{1 - z_y^2}}{\sin z}
	\]		
and the wave equation $z_{xy} = 0$ admit no rank-$1$ B\"acklund transformation relating
their solutions, but a rank-$2$ B\"acklund transformation does exist.		
		
In a future paper, we will present a partial classification of homogeneous rank-$2$ B\"acklund transformations relating two hyperbolic Monge-Amp\`ere systems. (It will turn out that the rank-$2$ B\"acklund transformation
corresponding to \eqref{Tztrans} is nonhomogeneous.)	Based on our classification so far, we expect that those
homogeneous B\"acklund transformations (relating two hyperbolic Monge-Amp\`ere systems) 
		that are `genuinely' rank-$2$
		are quite few. 		
		
\section{Acknowledgement}			
		
The author would like to thank his PhD thesis advisor, Prof. Robert L. Bryant, for all his guidance and support.
He would like to thank Prof. Jeanne N. Clelland for her advice on the current work.
Thanks to the referee for reading the manuscript and providing helpful suggestions and comments.

\begin{appendix}
\section{Calculations in Theorem \ref{generalitytheorem}}\label{App:genThm}
					
		This Appendix supplements the proof of Theorem \ref{generalitytheorem} by providing more 
		details of calculation. Most calculations below are computed using Maple\texttrademark.		
		
		First consider the case when, on $U$, $\epsilon_1 = \epsilon_2 = 1$.
		Since $P_{24}, P_{33},$ $P_{66}, P_{75}$ never appear in the
		 equation \eqref{rank1streq}, we can set them all to zero.
	 	Since $P_{14}$ and $P_{23}$ only appear in the term $(P_{14} - P_{23})\omega^3\W\omega^4$, 
		we can set $P_{14} = 0$. For similar reasons, we can set $P_{13}, P_{55}, P_{56} = 0$. For convenience, we
			rename $A_1$ as $P_{81}$ and $A_4$ as $P_{84}$.
		Now there are $42$ functions $P_{ij}$ remaining, and they are determined. 
			
		For each $P_{ij}$, there exist functions $P_{ijk}$ defined on $U$ satisfying
							\[
								\ed(P_{ij}) =  P_{ijk}\omega^k.
							\]
		We call these $P_{ijk}$ the \emph{derivatives}
		of $P_{ij}$.
		
		Now, applying $\ed^2 = 0$ to the equation \eqref{rank1streq}, we obtain
		$106$ polynomial equations expressed in terms of
		all $42$ $P_{ij}$ and $186$ of all $252$ $P_{ijk}$.
		These equations imply:
					\begin{equation*}
						\begin{alignedat}{3}
						P_{01} &= P_{41} - P_{51}, &&{\quad} P_{02} = P_{42} - P_{52}, &&{\quad} P_{03} =P_{53} -P_{43} - P_{81}, \\
						 P_{04} &= P_{54}-P_{44},&&{\quad}
						P_{05} =P_{15} -P_{45} - P_{84}, &&{\quad} P_{11}= -P_{51},\\
						 P_{12} &= -P_{52}, &&{\quad} P_{21} = P_{84}, &&{\quad} P_{22} = -1,\\
						 P_{35} &= 0, &&{\quad} P_{36} = -1, &&{\quad} P_{61} = 1, \\
						  P_{62} &= -P_{81}, &&{\quad} P_{73} = 0, &&{\quad} P_{74} = -1.
						  \end{alignedat}
					\end{equation*}
					
			With these relations, all coefficients in \eqref{rank1streq} can be expressed
			in terms of $27$ $P_{ij}$. 
			Repeating the steps above by defining the derivatives $P_{ijk}$ (now $162$ in all) 
			and applying $\ed^2 =0$ to \eqref{rank1streq}, 
			we obtain a system of $91$ polynomial equations in these $27$ $P_{ij}$ and
			$124$ of all 162 $P_{ijk}$, which imply
					\begin{align*}
						P_{31} &= -P_{32}P_{84} - P_{15} - P_{34} - 2P_{43} - 2P_{45}+P_{53}+P_{76} - P_{81} - P_{84}, \\
						P_{72}& = -P_{71}P_{81} - P_{15} +P_{34} +2P_{43} +2P_{45}-3P_{53} - P_{76}+P_{81}+P_{84}.
					\end{align*}
			
			Using these relations and repeating the steps above, we obtain
					\[
						P_{06} = P_{16} - P_{46}.
					\]		
			 All coefficients in \eqref{rank1streq} are then expressed in terms of $24$ $P_{ij}$.
			
			Now, corresponding to the remaining $24$ $P_{ij}$ are $144$ derivatives $P_{ijk}$. 
			Applying $\ed^2 = 0$
			to \eqref{rank1streq}
			yields a system of $88$ polynomial equations, expressed in terms of the $24$ $P_{ij}$ and $122$ of
			the $144$ derivatives $P_{ijk}$. This system can be solved for $P_{ijk}$; in the solution, all $P_{ijk}$ are
			expressed explicitly in terms of the $24$ $P_{ij}$ and $64$ $P_{ijk}$ that are `free'.

			Let ${a} = (a^\alpha)$ $(\alpha = 1,\ldots,24)$ stand for the $24$ remaining $P_{ij}$; let ${b} = (b^\rho)$ $(\rho = 1,\ldots,64)$ stand for the $64$ `free' $P_{ijk}$. We already have
					\begin{align}
						\ed\omega^i & = -\frac{1}{2}C^i_{jk}(a)\omega^j\W\omega^k, \label{estreqn0}\\
						\ed a^\alpha &= F^\alpha_i({a},{b}) \omega^i,   \label{estreqn}
					\end{align}
			for some real analytic functions $F^\alpha_i$ and $C^i_{jk}$ satisfying $C^i_{jk}+C^i_{kj} = 0$. 
			
			Now compute the exterior derivatives
					\[
						\ed(F^\alpha_i(a,b)\omega^i), \quad \alpha = 1,\ldots,24,
					\]	
			and take into account \eqref{estreqn0} and \eqref{estreqn}. From this we obtain $2$-forms
			$\Omega^\alpha$ 
			that are linear combinations of $\ed b^\rho \W\omega^i$ and $\omega^i\W\omega^j$.
			Let $\hat\Omega^\alpha$ denote the part of $\Omega^\alpha$ consisting of 
			linear combinations of $\ed b^\rho\W\omega^i$ only.
			Replacing $\ed b^\rho$ in $\hat\Omega^\alpha$ by $G^\rho_i\omega^i$ defines a linear map
					\[
						\phi: \Hom(\R^6,\R^{64})\rightarrow \Lambda^2(\R^6)^*\otimes \R^{24}
					\]
			at each point of $U$.
			
			Let $[\Omega]$ denote the equivalence class of $(\Omega^\alpha)$ in the cokernel of $\phi$.
			One can show that $[\Omega]$ must vanish and that its vanishing leads to a system of
			$35$ equations
			for $a$ and $b$.
			This system can be solved for $12$ of the $64$ components of $b$. Apply such a solution
			and update $a^\alpha$, $b^\rho$ and
			the functions $F^\alpha_i$ accordingly.
			
			It is not difficult to verify, using Maple\texttrademark, 
			that the updated $a^\alpha$ $(\alpha = 1,\ldots,24)$, $b^\rho$ $(\rho
			  = 1,\ldots,52)$, $C^i_{jk}$ and $F^\alpha_i$ satisfy the conditions {\bf(A)}-{\bf(C)}
			  in {\bf Step 1}.

			For {\bf Steps 2 and 3}, calculation shows that the \emph{tableaux of free derivatives}
			has Cartan
			characters 
					\[
						(s_1,s_2,s_3,s_4,s_5,s_6) = (24,22,6,0,0,0)
					\]
			and the dimension of its first prolongation 
					\[
						\delta = 64< s_1+2s_2+3s_3+4s_4+5s_5+6s_6 = 86.
					\]

			The cases when, on $U$, $\epsilon_1$ and $\epsilon_2$ take other values
			follow similar steps. In each of these cases,
			the last nonzero Cartan character, computed at a corresponding stage, is $s_3 = 6$.
			
	\section{Invariants of an Euler-Lagrange System}
\label{App:ELSym}
					
	This Appendix supplements the proof of Proposition \ref{Cohomog1SymmetryProp}.

	We start with the $G_1$-structure $\pi:\mathcal{G}_1\rightarrow M$ 
	of a hyperbolic Monge-Amp\`ere system~$(M,\I)$ (see \cite{BGG} or 
	Section \ref{InvMA}). Assume that $S_2 =\bs0$ (i.e., the Euler-Lagrange case).
	
	Recall that  the $2\times 2$-matrix $S_1: \mathcal{G}_1\rightarrow \gl(2,\R)$ is equivariant
	under the $G_1$-action. By \eqref{S1S2trans} and \eqref{S1S2trans2}, it is easy to see that, when
	$\det(S_1(u)) > 0$ (resp., $\det(S_1(u))<0$) at $u\in \mathcal{G}_1$, the same is true for
	$\det(S_1(u\cdot g))$ for all $g\in G_1$, and the matrix 
	$S_1(u)$ lies in the same $G_1$-orbit as $\diag(1,1)$
	(resp., $\diag(1,-1)$). 
	
	Now assume that $\det(S_1)>0$ holds on $\pi^{-1}U\subset \mathcal{G}$ for some domain
	 $U\subset M$. By the discussion above, we can reduce to a
	subbundle $\mathcal{H}\subset \mathcal{G}_1$ defined by $S_1 = \diag(1,1)$.
	
	It is easy to see that $\mathcal{H}$ is an $H$-structure on $U$ where
	\[
		H = \left\{\left.\left(\begin{array}{ccc}
						\epsilon&0&0\\
						0&A&0\\
						0&0&\epsilon A
			\end{array}\right)\right|\epsilon = \pm 1, ~A\in \GL(2,\R), ~\det(A) = \epsilon\right\}\subset G_1
	\]
	is a (disconnected) $3$-dimensional Lie subgroup.
	Let the restriction of $\pi: \mathcal{G}_1\rightarrow M$ 
	to $\mathcal{H}$ be denoted
	by the same symbol $\pi$. 
	
	One can verify that, restricted to $\mathcal{H}$, 		
		the $1$-forms $\phi_7-\phi_3$, $\phi_6-\phi_2$, 
		$\phi_5-\phi_1$ and $\phi_0$ in the equation \eqref{StrEqnMA}
		become semi-basic relative to $\pi:\mathcal{H}\rightarrow U$.
		Hence, there exist functions $Q_{ij}$ defined on $\mathcal{H}$ such that
			\begin{equation}
				\begin{array}{ll}
				\phi_7 = \phi_3+Q_{7i}\omega^i, \quad &\phi_6 = \phi_2+Q_{6i}\omega^i, \\[0.8em]
				\phi_5  = \phi_1+Q_{5i}\omega^i,\quad &\phi_0 = Q_{0i}\omega^i,
				\end{array}\label{betweenphis}
			\end{equation} 
		where the summations are over $i = 0,1,\ldots,4$.
		There are ambiguities in these $Q_{ij}$ as we can modify them without changing the form
		of the structure equation \eqref{StrEqnMA}. Using such ambiguities, we can arrange that
		 \begin{equation}
		 	Q_{71} =  Q_{73} = Q_{62} =  Q_{64} =  Q_{51} = Q_{52} =  Q_{53} = Q_{54} = 0;\label{NormalizeQij}
		\end{equation}	
		 the remaining $Q_{ij}$ are then determined.
		
		Applying $\ed^2 = 0$ to \eqref{StrEqnMA} and reducing appropriately, we obtain
			\begin{align*}
				\ed^2\omega^1 &\equiv{\phantom{-}} (Q_{63} - Q_{04})\omega^0\W\omega^3\W\omega^4 \mod \omega^1,\omega^2,\\
				\ed^2\omega^2 &\equiv {\phantom{-}} (Q_{03}-Q_{74})\omega^0\W\omega^3\W\omega^4 \mod
\omega^1,\omega^2,\\
				\ed^2\omega^3& \equiv {\phantom{-}} (Q_{02}+Q_{61})\omega^0\W\omega^1\W\omega^2\mod\omega^3,\omega^4,\\
				\ed^2\omega^4&\equiv (-Q_{01} - Q_{72})\omega^0\W\omega^1\W\omega^2\mod\omega^3,\omega^4.
			\end{align*}
		This implies that
				\[
					Q_{61} = -Q_{02}, \quad Q_{63} = Q_{04}, \quad Q_{72} = -Q_{01}, \quad Q_{74} = Q_{03}.
				\]
		Now all coefficients in the structure equation \eqref{StrEqnMA} 
		are expressed in terms of $Q_{0i}$ $(i = 0,1,\ldots,4)$ and $Q_{j0}$ $(j = 5,6,7)$.
		By applying $\ed^2 = 0$ to \eqref{StrEqnMA}, it is not difficult to verify that, reduced modulo
		$\omega^0,\omega^1,\ldots,\omega^4$, the following congruences hold:
			\begin{equation}\label{Qijtrans}
				\begin{alignedat}{1}
				\ed\left(\begin{array}{cc}
						Q_{01} &Q_{03}\\
						Q_{02} &Q_{04}
					\end{array}\right)&\equiv \left(\begin{array}{cc}
												\phi_1&\phi_3\\
												\phi_2&-\phi_1
											\end{array}
											\right)
											\left(\begin{array}{cc}
														Q_{01}&Q_{03}\\
														Q_{02}&Q_{04}
										\end{array}\right),\quad \ed Q_{00} \equiv 0,\\
				\ed\left(\begin{array}{ccc}
						Q_{50}\\
						Q_{60}\\
						Q_{70}
					\end{array}\right)&\equiv \left(\begin{array}{ccc}
												0& \phi_3&-\phi_2\\
												2\phi_2&-2\phi_1&0\\
												-2\phi_3&0&2\phi_1
											\end{array}\right)	
											\left(\begin{array}{ccc}
													Q_{50}\\
													Q_{60}\\
													Q_{70}
					\end{array}\right).
					\end{alignedat}
			\end{equation}	
		
		The congruences \eqref{Qijtrans} 
		tell us how the remaining $Q_{ij}$
		transform under the action by the identity component of $H$. Moreover, 
		it is easy to compute directly from \eqref{StrEqnMA} to verify that
			\begin{equation}\label{Qijdiscretetrans}
				\begin{alignedat}{1}
				Q_{00}(u\cdot h_0) &= -Q_{00}(u),\\
				\left(\begin{array}{cc}
						Q_{01} &Q_{03}\\
						Q_{02} &Q_{04}
					\end{array}\right)(u\cdot h_0)& = \left(\begin{array}{cc}
						-Q_{01} &Q_{03}\\
						Q_{02} &-Q_{04}
					\end{array}\right)(u),  \\
			\left(\begin{array}{c}Q_{50}\\ Q_{60}\\ Q_{70}\end{array}\right)(u\cdot h_0) 
			&= \left(\begin{array}{c}-Q_{50}\\ Q_{60}\\ Q_{70}\end{array}\right)(u), \\[0.5em]
			   h_0 = \diag(-1,&-1,1,1,-1)\in H	
			  	\end{alignedat}
			\end{equation}
		hold for any $u\in \mathcal{H}$. 
		
		Note that ${H}$ is generated by its identity component
		and $h_0$.
		Combining \eqref{Qijtrans} and \eqref{Qijdiscretetrans},
		it is easy to see that $Q_{01}Q_{04} - Q_{02}Q_{03}$ and $|Q_{00}|$ are local invariants  
		of the underlying Euler-Lagrange system. 
		
		Moreover, using \eqref{Qijtrans} and \eqref{Qijdiscretetrans}, it is easy to see that the 
		$H$-orbit of 
		\[
		q(u) :=\left(\begin{array}{cc}
														Q_{01}&Q_{03}\\
														Q_{02}&Q_{04}
												\end{array}\right)(u), \quad u\in \mathcal{H}
		\]										
		consists of all $2$-by-$2$ matrices with the same determinant as $q(u)$. Now we are ready to prove
		the following lemma.
		
	{\lemma If $\det(q)\ne 0$ on $\mathcal{H}$,
			then there is a canonical way to define a coframing on $U$. \label{AppLemmaCoframing}}\vskip 2mm
		
		{\it Proof}. If the function $L:=\det(q)$ is nonvanishing on $U$, 
		one can reduce to the subbundle $\mathcal{H}_1$ of $\mathcal{H}$ defined by 
		 $q = \diag(L,1)$.
		It is easy to see that each fiber of $\mathcal{H}_1$ over $U$ contains a single element.
		\qed
		
	{\remark	As a result of Lemma \ref{AppLemmaCoframing}, if $\det(q)\ne 0$ on $U$, 
		then the
		corresponding hyperbolic Euler-Lagrange system has a symmetry of dimension at most $5$.
		This is a consequence of applying the \emph{Frobenius Theorem}.
		\label{rmksymmetrydimension}}\\

	Now we proceed to complete the proof of Proposition \ref{Cohomog1SymmetryProp}.
	Recall that the coframing
	$(\eta^0,\eta^1,\ldots,\eta^4)$ and the $\phi_\alpha$ in \eqref{etasandphis} 
	verify the equation \eqref{StrEqnMA},
	$S_1 = \diag(1,1)$, and $S_2 = \bs 0$. Moreover, we have chosen 
	the $\phi_\alpha$ to satisfy
	\eqref{NormalizeQij}, where $Q_{ij}$ are computed using \eqref{betweenphis}. 
	By \eqref{etasandphis}, it is immediate that
			\begin{align*}
				Q_{00} = Q_{02}& =0, \quad Q_{01} = -Q_{04} =\frac{1}{\sqrt{2}}, \quad Q_{03} = 1,\\
					Q_{70}& = 1, \quad Q_{60}=-1, \quad Q_{50} = \sqrt{2}.
			\end{align*}
		Clearly, $\det(q) = Q_{01}Q_{04} - Q_{02}Q_{03}= -1/2\ne 0$. 
		By Lemma \ref{AppLemmaCoframing} and Remark \ref{rmksymmetrydimension}, the 
		hyperbolic Euler-Lagrange system considered in Proposition  \ref{Cohomog1SymmetryProp} has a 
		symmetry of dimension at most $5$. Because that Euler-Lagrange system is homogeneous,
		it follows that its symmetry has dimension $5$.

\end{appendix}

\bibliographystyle{alpha}
\normalbaselines 
\newcommand{\etalchar}[1]{$^{#1}$}

\end{document}